\newcommand{\papertitle}{Learning Convex Optimization Control Policies}

\documentclass[12pt]{article}
\usepackage{fullpage}

\title{\papertitle{}}
\author{Akshay Agrawal \and
Shane Barratt \and
Stephen Boyd \and
Bartolomeo Stellato\thanks{Authors listed in alphabetical order.}}

\newcommand{\citep}{\cite}

\usepackage[colorlinks,allcolors=blue,bookmarks=false,hypertexnames=true]{hyperref}
\usepackage{adjustbox}
\usepackage{mathtools,graphicx,wrapfig,graphics,psfrag,amsmath,amsfonts,verbatim,xcolor,color,siunitx}

\newcommand{\BEAS}{\begin{eqnarray}}
\newcommand{\EEAS}{\end{eqnarray}}
\newcommand{\BEA}{\begin{eqnarray}}
\newcommand{\EEA}{\end{eqnarray}}
\newcommand{\BEQ}{\begin{equation}}
\newcommand{\EEQ}{\end{equation}}
\newcommand{\BIT}{\begin{itemize}}
\newcommand{\EIT}{\end{itemize}}
\newcommand{\BNUM}{\begin{enumerate}}
\newcommand{\ENUM}{\end{enumerate}}

\newcommand{\BA}{\begin{array}}
\newcommand{\EA}{\end{array}}

\newcommand{\eg}{{\it e.g.}}
\newcommand{\ie}{{\it i.e.}}

\newcommand{\ones}{\mathbf 1}

\newcommand{\reals}{{\mbox{\bf R}}}

\newcommand{\symm}{{\mbox{\bf S}}}  %

\newcommand{\Expect}{\mathop{\bf E{}}}

\newcommand{\argmin}{\mathop{\rm argmin}}

\newcommand{\var}{\mathop{\bf var}}

\newcommand{\argmax}{\mathop{\rm argmax}}

\newcounter{algorithmctr}[section]
\renewcommand{\thealgorithmctr}{\thesection.\arabic{algorithmctr}}
   {\refstepcounter{algorithmctr}\begin{list}{}{%
       \setlength{\rightmargin}{0.03\linewidth}%
       \setlength{\leftmargin}{0.03\linewidth}}%
       \rmfamily\small
       \item[]{\setlength{\parskip}{0ex}\hrulefill\par%
        \nopagebreak{\bfseries\textsf{Algorithm \thealgorithmctr~}}}}%
   {{\setlength{\parskip}{-3ex}\nopagebreak\par\hrulefill} \end{list}}

\begin{document}
\maketitle

\begin{abstract}
Many control policies used in various applications determine the input or
action by solving a convex optimization problem that depends on the current state
and some parameters. Common examples of such
convex optimization control policies (COCPs) include the linear quadratic regulator (LQR),
convex model predictive control (MPC), and convex control-Lyapunov or approximate
dynamic programming (ADP) policies.
These types of control policies are tuned by varying the parameters in the
optimization problem, such as the
LQR weights, to obtain good performance, judged by application-specific metrics.
Tuning is often done by hand, or by simple methods such as a crude grid search.
In this paper we propose a method to automate this process, by adjusting the
parameters using an approximate gradient of the performance metric with respect
to the parameters. Our method relies on recently developed methods that can
efficiently evaluate the derivative of the solution of a convex optimization
problem with respect to its parameters. We illustrate our method on several
examples.
\end{abstract}

\section{Introduction}
\subsection{Convex optimization control policies}

We consider the control of a stochastic dynamical system with known dynamics,
using a control policy that determines the input or action by
solving a convex optimization problem.  We call such policies
\emph{convex optimization control policies} (COCPs).
Many practical policies have this form, including the first modern
control policy,
the linear quadratic regulator~(LQR)~\citep{kalman1960contributions}.
In LQR, the convex optimization problem has quadratic objective and linear equality
constraints, and so can be solved explicitly, yielding the familiar
linear control policy.
More modern examples, which rely on more complicated optimization problems such as
quadratic programs (QPs),
include convex model predictive control (MPC)~\citep{borrelli2017predictive}
and convex approximate dynamic programming (ADP)~\citep{bertsekas2017dynamic}.
These policies are used in many applications, including
robotics~\citep{Kuindersma14},
vehicle control~\citep{stewart2008},
rocket landing~\citep{spacex},
supply chain optimization~\citep{powell2012logistics},
and finance~\citep{markowitz1952portfolio,cornuejols2006,boyd2017multi}.

Control policies in general, and COCPs in particular, are judged by
application-specific metrics; these metrics are evaluated using simulation with
historical or simulated values of the unknown quantities.  In some but not all
cases, the metrics have the traditional form of the average value of a given
stage cost. We consider here more general metrics that can be functions of the
whole state and input trajectories. An example of such a metric is the expected
drawdown of a portfolio over some time period, \ie, the expected value of the
minimum
future value of a portfolio.

In a few cases, the optimal policy for a traditional stochastic control problem
has COCP form.  A well-known example is LQR~\citep{kalman1960contributions}.
Another generic example is when the dynamics are affine
and the stage cost is convex, in which case the Bellman value function is convex,
and evaluating the optimal policy reduces to solving a convex optimization problem~\citep[\S3.3.1]{keshavarz2012convex}.
While it is nice to know that in this case that the optimal policy has COCP form,
we generally cannot express the value function in a form
that allows us to evaluate the policy, so this observation is not useful
in practice.
In a far wider set of cases, a COCP policy is not optimal, but only a
good, practical heuristic.

COCPs have some attractive properties compared to other parametrized
control policies.  When the convex problem to be solved is well chosen,
the policy is at least reasonable for any choice of the parameter values
over the allowed set.  As a specific example, consider a linear control policy
parametrized by the gain matrix, which indeed would seem to be the most
natural parametrization of a linear policy.
The set of gain matrices that lead to a stable closed-loop system
(a very minimal performance requirement) can be very complex, even disconnected.
In contrast, consider an LQR control policy parametrized by a state and control
cost matrix (constrained to be positive definite).  In this case any choice of policy
yields a stable closed-loop system.
It is far easier and safer to tune parameters when any feasible choice leads to
at least a reasonable policy.

All control policies are tuned by choosing various parameters that appear
in them.  In the case of COCPs, the parameters are in the optimization problem
that is solved to evaluate the policy.
The tuning is usually done based on simulation with historical disturbances
(called \emph{back-testing}) or synthetic disturbances.
It is often done by hand, or by a crude grid search.
A familiar example of this is tuning the weights in an LQR controller to
obtain good practical performance~\citep{anderson1990}.

In this paper we present an automated method for tuning parameters
in COCPs to achieve good values of a performance metric.
Our method simulates the closed-loop system, \ie, the system with the
policy in the loop,
and computes an approximate (stochastic) gradient of the expected performance
with respect to the parameters.  It uses this gradient to update the parameters
via a projected stochastic gradient method.
Central to our method is the
fact that the solution map for convex optimization problems is often
differentiable, and its derivative can be efficiently computed
\citep{agrawal2019differentiable,amos2019differentiable}.
This is combined with relatively new implementations of automatic differentiation,
widely used in training neural networks \citep{abadi2016tensorflow, paszke2019pytorch}.

Our method is not guaranteed to find the best parameter values, since the
performance metric is not a convex function of the COCP parameter values, and
we use a local search method. This is not a problem in practice, since in a
typical use case, the COCP is initialized with reasonable parameters,
and our method is used to tune these parameters to improve the performance
(sometimes considerably).

\subsection{Related work}

\paragraph{Dynamic programming.}
The Markov decision process (MDP)
is a general stochastic control problem that can be solved in principle
using dynamic programming (DP) \citep{Bellman:1957,bellman1957markovian,bertsekas2017dynamic}.
The optimal policy is evaluated by solving an optimization problem,
one that includes a current stage cost and the expected value of cost-to-go
or value function at the next state.
This optimization problem corresponds to a COCP when the system dynamics are 
linear or affine and the stage cost is convex~\cite{bertsekas2017dynamic}.
Unfortunately, the value function can be found in a tractable form in only
a few cases.
A notable tractable case is when the cost is a convex extended
quadratic and the dynamics are affine~\citep{barratt2018stochastic}.

\paragraph{Approximate dynamic programming.}
ADP~\citep{bertsekas1996neuro,powell2007approximate} refers to heuristic
solution methods for stochastic control problems that
replace the value function in DP with an approximation, or
search over a parametric family of policies~\citep[\S2.1]{bertsekas2019reinforcement}.

In many ADP methods, an offline
optimization problem is solved to approximate the value function. When there are a
finite number of state and inputs, the approximation problem can be written as
a linear program (LP) by relaxing the Bellman equation to an
inequality~\citep{de2003linear}. When the dynamics are linear, the cost is
quadratic, and the input is constrained to lie in a convex set,
an approximate convex quadratic value function can be found
by solving a particular semidefinite program (SDP)~\citep{wang2009performance}.
The quality of the approximation can
also be improved by iterating the Bellman
inequality~\citep{wang2015approximate, stellato2017}. Because the approximate
value function is convex quadratic and the dynamics are linear, the resulting policy is a COCP.

Other methods approximate the cost-to-go by iteratively adjusting the approximate value function
to satisfy the Bellman equation.
Examples of these methods include projected value iteration or fitted
Q-iteration~\citep{gordon1995stable}, temporal difference
learning~\citep{sutton1988learning,bertsekas2004improved}, and approximate
policy iteration~\citep{nedic2003least}.
Notable applications of COCPs here include the use of quadratic approximate
cost-to-go functions for input-affine systems with convex cost, which can be approximately fit
using projected value iteration~\citep{keshavarz2014quadratic}, and modeling the state-action cost-to-go
function as an input-convex neural network~\cite[\S6.4]{amos2017input}. Other
approximation schemes fit nonconvex value functions, so the resulting policies
are not necessarily COCPs.
Notably, when the parametrization involves a
featurization computed by a deep neural network, the ADP method is an instance
of deep reinforcement learning.

Other ADP methods parametrize the policy and tune the parameters directly
to improve performance; this is often referred to as policy search or policy
approximation \citep[\S5.7]{bertsekas2019reinforcement}.
The most common method is gradient or stochastic gradient search
\citep[\S7.2]{powell2007approximate}, which is the method we
employ in this paper, with a parametrized COCP as the policy.
Historically, the
most widely used of these policy approximation methods is the
Proportional-Integral-Derivative (PID) controller
\citep{minorsky1922directional}, which indeed
can be tuned using gradient methods \citep{aastrom1993automatic}.

\paragraph{Reinforcement learning.}
Reinforcement learning (RL)~\citep{sutton2018reinforcement} and adaptive
control~\citep{aastrom2013adaptive} are essentially equivalent to ADP
\citep[\S1.4]{bertsekas2019reinforcement}, but with different notation and
different emphasis. RL pays special attention to problems in
which one does not possess a mathematical model of the dynamics or the
expected cost, but has access to a computational simulator for both.
Our method cannot be used directly in this setting, since we assume that we
have mathematical descriptions of the dynamics and cost. However, our method
might be used after learning a suitable model of the dynamics and cost. Alternatively,
COCPs could be used as part of the policy in modern policy gradient or
actor-critic
algorithms~\citep{williams1987reinforcement,lillicrap2015continuous,schulman2017proximal}.

\paragraph{Learning optimization-based policies.} Other work
has considered tuning optimization-based control policies.
For example, there is prior
work on learning for MPC, including
nonconvex MPC controllers \citep{amos2018mpc}, cost function shaping
\citep{tamar2017learning}, differentiable path integral control
\citep{okada2017path}, and system identification of terminal constraint sets
and costs \citep{rosolia2017learning}. As far as we are aware, our
work is the first to consider the specific class of parametrized convex
programs.

\paragraph{Real-time optimization.}
COCPs might be considered computationally expensive control policies compared
to conventional analytical control policies such as the linear control policy
prescribed by LQR. However, this is not the case in practice, thanks to fast
embedded solvers~\citep{domahidi2013ecos,osqp,wang2010fast} and code generation
tools that emit solvers specialized to parametric problems
~\citep{Mattingley:2012,chu2013code,osqp_codegen}. For example, the aerospace
and space transportation company SpaceX uses the QP code generation tool
CVXGEN~\cite{Mattingley:2012} to land its rockets~\citep{spacex}. COCPs based
on MPC, which have many more variables and constraints than those based on ADP,
can also be evaluated very efficiently~\citep{bemporad2002,wang2010}, even at
MHz rates~\citep{jerez2014}.

\subsection{Outline}
In \S\ref{sec:prob}, we introduce
the controller tuning problem that
we wish to solve.
In \S\ref{sec:examples_cocps}, we describe
some common forms of COCPs.
In \S\ref{sec:solution_method}, we propose
a heuristic for the controller tuning problem.
In \S\ref{sec:examples}, we apply our heuristic for tuning COCPs
to examples in portfolio optimization, vehicle control,
and supply-chain management. We conclude in
\S\ref{sec:extensions_and_variations} by discussing extensions
and variations.

\section{Controller tuning problem}\label{sec:prob}

\paragraph{System dynamics.}
We consider a dynamical system with dynamics given by
\BEQ
x_{t+1} = f(x_t, u_t, w_t), \quad t=0,1, \ldots.
\label{eq:dyn}
\EEQ
At time period $t$, $x_t\in\reals^n$ is the state,
$u_t \in \reals^m$ is the input or action,
$w_t\in\mathcal W$ is the disturbance,
and $f:\reals^n\times\reals^m\times\mathcal W\to\reals^n$ is
the state transition function.
The initial state $x_0$ and the disturbances $w_t$ are random variables.
In the traditional stochastic control problem, it is
assumed that $x_0, w_0, w_1, \ldots$ are independent, with
$w_0, w_1, \ldots $ identically
distributed.  We do not make this assumption.

The inputs are given by a state feedback control policy,
\BEQ\label{e-policy}
u_t=\phi(x_t),
\EEQ
where $\phi : \reals^n \to \reals^m$
is the policy.
In particular, we assume the state $x_t$ at time period $t$ is fully observable
when the input $u_t$ is chosen.
It will be clear later that this assumption is not really needed, since
our method can be applied to an estimated state feedback policy,
either with a fixed state estimator, or with a state estimator that is
itself a parametrized convex problem (see \S\ref{sec:extensions_and_variations}).

With the dynamics~(\ref{eq:dyn}) and policy~(\ref{e-policy}), the state
and input trajectories $x_0, x_1, \ldots$ and $u_0, u_1, \ldots $ form
a stochastic process.

\paragraph{Convex optimization control policies.}
We specifically consider COCPs, which have the form
\BEQ\label{eq:cocp}
\begin{array}{llll}
\phi(x) &= & \underset{u}{\argmin} & f_0(x, u; \theta)  \\
&& \mbox{subject to} & f_i(x, u; \theta) \leq 0, \quad i=1,\ldots,k,\\
&& & g_i(x, u; \theta) = 0, \quad i=1,\ldots,\ell,
\end{array}
\EEQ
where $f_i$ are convex in $u$ and $g_i$ are affine in $u$. To evaluate a COCP
we must solve a convex optimization problem, which we assume has a unique
solution. The convex optimization problem \eqref{eq:cocp} is given by a
\emph{parametrized problem description} \citep[\S4.1.4]{boyd2004convex}, in
which the vector $\theta \in \Theta \subseteq \reals^p$ is the parameter ($\Theta$ is the set of
allowable parameter values). The value of the parameter $\theta$ (and $x$)
specifies a particular problem instance, and it can be adjusted to tune the
control policy. The problem we address in this paper is the choice of the
parameter $\theta$.

\paragraph{Performance metric.}
We judge the performance of a control policy, or choice of control
policy parameter $\theta$, by the average value of a cost over
trajectories of length $T$.
Here the horizon $T$ is chosen large enough so that the average over $T$ time steps is
close enough to the long term average.
We denote the trajectories over $t=0,\ldots, T$ as
\[
\begin{array}{rcl}
X &=& (x_0, x_1, \ldots, x_T) \in \reals^{N},\\
U &=& (u_0, u_1, \ldots, u_T)\in \reals^{M},\\
W &=& (w_0, w_1, \ldots, w_T)\in \mathcal W^{T+1},
\end{array}
\]
where $N=(T+1)n$ and $M=(T+1)m$.
These state, input, and disturbance trajectories are random variables, with
distributions that depend on the parameter $\theta$.

The cost is provided by a function
$\psi : \reals^N \times \reals^M \times \mathcal W^{T+1} \to \reals \cup \{+\infty\}$.
Infinite values of $\psi$ can be interpreted as encoding constraints on the trajectories.
A policy is judged by the expected value of this cost,
\[
J(\theta) = \Expect \psi(X, U, W).
\]
We emphasize that $J$ depends on the control policy parameter $\theta$, since
$x_1, \ldots, x_T$ and $u_0, \ldots, u_T$ depend on $\theta$.

We mention that the traditional cost function is separable, with the form
\BEQ\label{e-stage-cost}
\psi(X,U,W) = \frac{1}{T+1} \sum_{t=0}^T g(x_t,u_t,w_t),
\EEQ
where $g:\reals^n \times \reals^m \times \mathcal W \to \reals \cup \{\infty\}$ is a
stage cost function.
However, we do not require a cost function that is separable across time.

\paragraph{Evaluating $J(\theta)$.} We generally cannot evaluate $J(\theta)$
exactly. Instead, assuming that we can sample the initial state and the disturbances,
we can compute a Monte Carlo approximation of it.  In the simplest version,
we generate $K$ independent trajectories
\[(X^1, U^1, W^1), \ldots, (X^K, U^K, W^K),\]
and form the approximation
\[
\hat J(\theta) = \frac{1}{K} \sum_{i=1}^{K}\psi(X^i, U^i, W^i).
\]
This computation requires carrying out $K$ simulations over $T$ time steps,
which involves solving $K(T+1)$ convex optimization problems to evaluate $u_t^i$,
$t=0, \ldots, T$, $i=1,\ldots K$.

Evidently, $\hat J(\theta)$ is an unbiased
approximation of $J(\theta)$, meaning
\[
\Expect \hat J(\theta) = J(\theta).
\]
The quality of this approximation increases as $K$ increases, since
\[
\var{ \hat J(\theta) } = \frac{\var\psi(X, U, W)}{K},
\]
where $\var$ denotes the variance; \ie, the variance goes to $0$ as $K$ gets large.
Of course more sophisticated methods can be used to approximately
evaluate $J(\theta)$, \eg, importance sampling (see \citep{cochran2007sampling}).

\paragraph{Controller tuning problem.}
The controller tuning problem has the form
\begin{equation}
\begin{array}{ll}
\mbox{minimize} & J(\theta) \\
\mbox{subject to} & \theta \in \Theta,
\end{array}
\label{eq:objective}
\end{equation}
with variable $\theta$. This is the problem we seek to solve in this paper.

\section{Examples of COCPs}
\label{sec:examples_cocps}
In this section we describe some common COCPs.

\paragraph{Optimal (dynamic programming) policy.}
In the traditional stochastic control setting, the cost function
is the average of stage costs computed by a function $g$, as in~\eqref{e-stage-cost}, and
$x_0,w_0,w_1, \ldots$ are independent.
Under some technical conditions, the optimal policy for $T\to \infty$,
\ie, the policy that minimizes $J$ over all possible state feedback policies,
and not just those of COCP form, has the form
\BEQ\label{e-dp}
\phi(x) = \argmin_u \Expect \left( g(x,u,w) + V( f(x,u,w)) \right),
\EEQ
where $V : \reals^n \to \reals$ is the optimal cost-to-go or Bellman
value function.
This form of the optimal policy is sometimes called the
dynamic programming (DP) form.
When $f$ is affine in $x$ and $u$, and $g$ is convex in $x$ and $u$,
it can be shown that the value function $V$ is convex \citep[\S3.3.1]{keshavarz2012convex}, so
the expression to be minimized above is convex in $u$,
and the optimal policy has COCP form (with no parameter $\theta$).

Unfortunately the optimal value function $V$ can be expressed in
tractable form in only a few special cases.
One well-known one is LQR
\citep{kalman1960contributions}, which has dynamics and stage cost
\BEQ
\label{eq:lqr}
f(x,u,w)=Ax+Bu+w, \qquad
g(x,u,w)=x^TQx+u^TRu,
\EEQ
with $A \in \reals^{n \times n}$,
$B \in \reals^{n \times m}$,
$Q \in \symm_+^n$ (the set of $n\times n$ symmetric positive
semidefinite matrices),
$R \in \symm_{++}^m$ (the set of symmetric positive
definite matrices),
and
$w \sim \mathcal N(0,\Sigma)$.
In this special case we can compute the value function,
which is a convex quadratic $V(x) = x^TPx$,
and the optimal policy has the form
\[
\phi(x) = \argmin_u \left(u^TRu + (Ax+Bu)^T P (Ax+Bu) \right)
= Kx,
\]
with
\[
K = -(R+B^TPB)^{-1}B^TPA.
\]

Note that we can consider the policy above as a COCP, if we
consider $P$ as our parameter $\theta$ (constrained to be
positive semidefinite).
Another option is to take
$P=\theta^T\theta$, where $\theta \in \reals^{n \times n}$,
so the COCP has objective
\[
f_0(x,u;\theta) =
u^TRu + \| \theta (Ax+Bu) \|_2^2.
\]

\paragraph{Approximate dynamic programming policy.}
An ADP \citep{powell2007approximate} or control-Lyapunov \citep{corless1988controller} policy has the form
\BEQ\label{e-adp}
\phi(x) = \argmin_u \Expect (g(x,u,w) + \hat V( f(x,u,w))),
\EEQ
where $\hat V$ is an approximation of the value function
for which the minimization over $u$ above is tractable.
When $g$ is convex in $u$, $f$ is affine in $u$, and $\hat V$ is
convex, the minimization above is a convex optimization problem~\cite{boyd2004convex}.
With a suitable parametrization of $\hat V$, this policy has
COCP form~\cite{keshavarz2012convex}.

\paragraph{Model predictive control policy.}
Suppose the cost function has the form~\eqref{e-stage-cost}, with stage
cost $g$.
In an MPC policy, the input is determined by solving an approximation
to the control problem over a short horizon, where the unknown
disturbances are replaced by predictions~\citep{rawlings2009model},
and applying only the first input.
A terminal cost function $g_H$ is often included in the optimization.

An MPC policy has the form
\begin{equation*}
\begin{array}{llll}
\phi(x) &=& \underset{u_0}{\argmin} & \sum_{t=0}^{H-1} g(x_t,u_t,\hat w_t) +
g_H(x_H) \\
&&\mbox{subject to} & x_{t+1} = f(x_t,u_t,\hat w_t), \quad t=0,\ldots,H-1,\\
&&& x_0 = x,
\end{array}
\end{equation*}
where $H$ is the planning horizon and
$\hat w_0,\ldots,\hat w_{H-1}$ are the predicted disturbances.
This optimization problem has variables $u_0,\ldots,u_{H-1}$
and $x_0,\ldots,x_H$; however, the $\argmin$ is over $u_0$ since in
MPC we only apply the first input.

When $f$ is affine in $(x,u)$, $g$ is convex in $(x,u)$,
and the terminal cost function $g_H$ is convex,
the minimization above is a convex optimization problem.
With a suitable parametrization of the terminal cost function $g_H$,
the MPC policy has COCP form.
When $f$ is not affine or $g$ is not convex, they can
be replaced with parametrized convex approximations.
The function that predicts the disturbances can also
be parametrized~(see~\S\ref{sec:extensions_and_variations}).

\section{Solution method}
\label{sec:solution_method}

Solving the controller tuning problem \eqref{eq:objective} exactly is in general
hard, especially when the number of parameters $p$ is large, so we will solve it approximately.
Historically, many practitioners have used derivative-free methods
to tune the parameters in control policies. Some of these methods include
CMA-ES~\citep{hansen2001completely} and other evolutionary strategies~\citep{salimans2017evolution}, Bayesian optimization~\citep{movckus1975bayesian}, grid search, and random search~\citep{anderson1953recent, solis1981minimization, bergstra2012random}. Many
more methods are catalogued in~\citep{conn2009introduction}. These methods can
certainly yield improvements over an initialization; however, they often
converge very slowly.

\paragraph{A gradient-based method.}
It is well-known that first-order optimization methods, which make use of
derivatives, can
outperform derivative-free methods. In this paper, we apply the projected
stochastic (sub)gradient method~\citep{robbins1951stochastic} to approximately
solve~\eqref{eq:objective}. That is, starting with initial parameters~$\theta^0$, at iteration~$k$, we simulate the system and compute~$\hat
J(\theta^k)$. We then compute an unbiased stochastic gradient of $J$,
$g^k = \nabla \hat J(\theta^k)$, by the
chain rule or backpropagation through time (BPTT)~\citep{rumelhart1988learning,werbos1990backpropagation}, and update the
parameters according to the rule $\theta^{k+1} =
\Pi_{\Theta}(\theta^k - \alpha^k g^k)$, where~$\Pi_\Theta(\theta)$
denotes the Euclidean projection of $\theta$ onto
$\Theta$ and $\alpha^k > 0$ is a step size.
Of course more sophisticated methods
can be used to update the parameters, for example,
those that employ momentum, variance reduction, or second-order information
(see~\citep{bottou2018optimization} and the references therein
for some of these methods).

\paragraph{Computing $g^k$.}
The computation of $g^k$ requires differentiating through the dynamics $f$,
the cost $\psi$, and, notably, the solution map $\phi$ of a convex optimization
problem. Methods for differentiating through special subclasses
of convex optimization have existed for many decades; for example, literature
on differentiating through QPs dates back to at least the
1960s \citep{boot1963sensitivity}. Similarly, it is well known that if
the objective function and constraint functions of a convex optimization
problem are all smooth, and some regularity conditions are satisfied, then its
derivative can be computed by differentiating through the KKT optimality
conditions \citep{kuhn1951nonlinear,barratt2018differentiability}. Until very recently, however, it was not
generically possible to differentiate through a convex optimization problem
with nondifferentiable objective or constraints; recent work~\citep{busseti2018solution, diffcp2019, agrawal2019differentiable,
amos2019differentiable} has shown how to efficiently
and easily compute this derivative.

\paragraph{Non-differentiability.}
Until this point, we have assumed the differentiability of all of the functions
involved ($f$, $\psi$, and $\phi$). In
real applications, these functions very well may not be differentiable
everywhere. So long as the functions are differentiable almost everywhere,
however, it is reasonable to speak of applying a projected stochastic gradient method
to~\eqref{eq:objective}. At non-differentiable points, we compute a
heuristic quantity. For example, at some non-differentiable points
of $\phi$, a certain matrix fails to be invertible, and we compute a
least-squares approximation of the derivative instead, as in \citep{diffcp2019}.
In this sense, we overload the notation $\nabla f(x)$ to denote a gradient when
$f$ is differentiable at $x$, or some heuristic quantity (a ``gradient'') when
$f$ is not differentiable at $x$. In practice, as our examples in
\S\ref{sec:examples} demonstrate, we find that this
method works well. Indeed, most neural networks that are trained today
are not differentiable (\eg, the rectified linear unit or positive
part is a nondifferentiable activation function that is widely used) or even
subdifferentiable (since neural networks are usually nonconvex), but it is nonetheless possible to
train them, successfully, using stochastic ``gradient'' descent~\cite{goodfellow2016deep}.

\section{Examples}\label{sec:examples}
In this section, we present examples that illustrate
our method. Our control policies were implemented using CVXPY
\citep{diamond2016cvxpy, agrawal2018rewriting}, and we used
cvxpylayers \citep{agrawal2019differentiable} and PyTorch
\citep{paszke2019pytorch} to differentiate through them; cvxpylayers uses the
open-source package SCS \citep{odonoghue2016conic, odonoghue2017scs} to solve
convex optimization problems. For each example, we give the dynamics, the cost,
the COCP under consideration, and the result of applying our method to a
numerical instance.

In the numerical instances, we pick the number of simulations $K$ so that
the variance of $\hat J(\theta)$ is sufficiently small, and we tune the
step-size schedule $\alpha^k$ for each problem. %
BPTT is susceptible to exploding and vanishing
gradients \citep{bengio1994learning}, which can make learning difficult. This
issue can be mitigated by gradient clipping and regularization
\citep{pascanu2013difficulty}, which we do in some of our experiments.

\subsection{LQR}
\label{sec:lqr_example}
We first apply our method to the classical LQR problem, with dynamics
and cost
\[
f(x,u,w) = Ax + Bu + w, \quad \psi(X, U, W) = \frac{1}{T+1}\sum_{t=0}^{T} x_t^TQx_t + u_t^TRu_t,
\]
where $A \in \reals^{n \times n}$,
$B \in \reals^{n \times m}$,
$Q \in \symm_+^n$,
$R \in \symm_{++}^m$,
and
$w \sim \mathcal N(0,\Sigma)$.

\paragraph{Policy.}
We use the COCP
\begin{equation} \label{eq:lqr_policy}
\phi(x) = \underset{u}\argmin \left(  u^TRu + \|\theta(Ax+Bu)\|_2^2 \right),
\end{equation}
with parameter $\theta \in\reals^{n\times n}$.  This policy is linear, of
the form $\phi(x) = Gx$, with
\[
G = -(R+B^T\theta^T\theta B)^{-1}B^T\theta^T\theta A.
\]
This COCP is clearly over-parametrized; for example, for any orthogonal
matrix $U$, $U\theta$ gives the identical policy as $\theta$.
If the matrix $\theta^T \theta$
satisfies a particular algebraic Riccati equation involving $A$, $B$, $Q$, and $R$, then
\eqref{eq:lqr_policy} is optimal (over all control policies) for the case $T\to \infty$.

\begin{figure}
  \centering
  \includegraphics[]{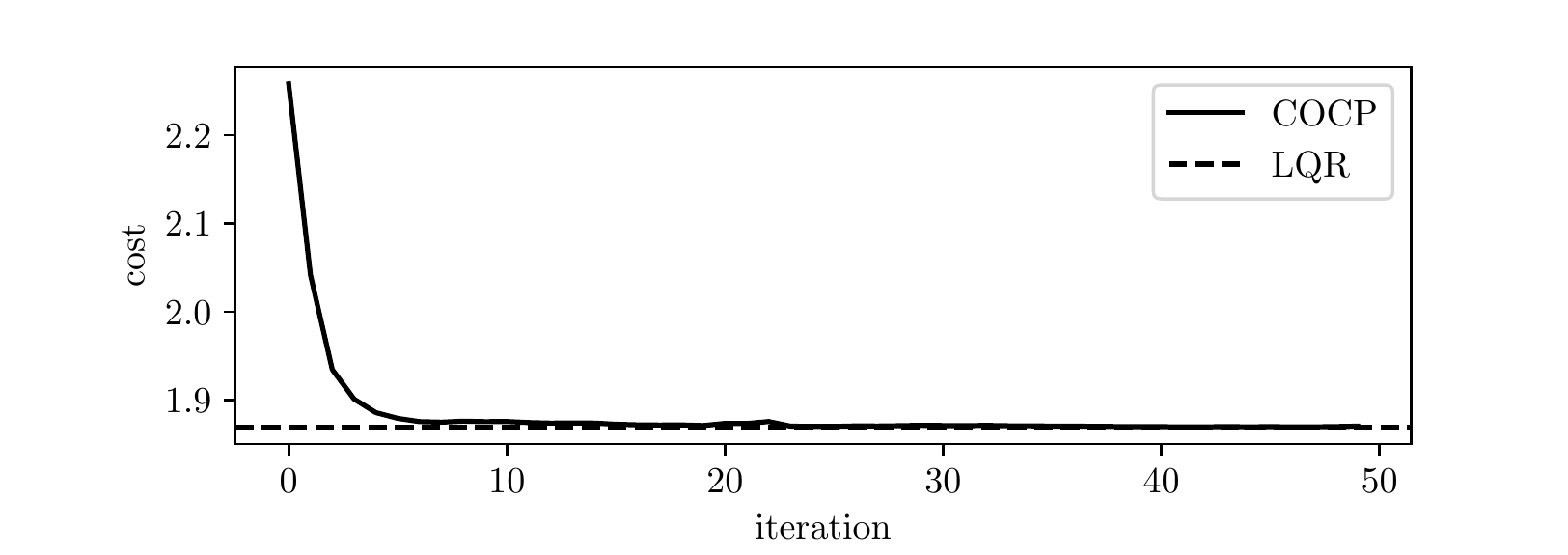}
  \caption{Tuning an LQR policy.}
  \label{fig:lqr}
\end{figure}

\paragraph{Numerical example.}
We consider a numerical example with
$n=4$ states, $m=2$ inputs, and $T=100$.
The entries of $A$ and $B$ were sampled from the standard normal distribution,
and we scaled $A$ such that its spectral radius was one. The cost matrices are
$Q=I$ and $R=I$, and the noise covariance is $W=(0.25)I$. We initialize~$\theta$ with the identity. We trained our policy~\eqref{eq:lqr_policy} for
50 iterations, using $K=6$ simulations per step,
starting with a step size of 0.5 that was decreased to 0.1 after 25 iterations.
Figure~\ref{fig:lqr} plots the average cost
of the COCP during learning versus the average cost of the
optimal LQR policy (in the case $T\to\infty$).
Our method appears to converge to near
the optimal cost in just 10 iterations.

\subsection{Box-constrained LQR}
A box-constrained LQR problem has the same dynamics and cost as LQR,
with an additional constraint $\|u_t\|_\infty \leq u_\mathrm{max}$:
\[
\psi(X, U, W) = \frac{1}{T+1}\sum_{t=0}^{T} g(x_t, u_t, w_t),  \qquad
g(x_t, u_t, w_t) = \begin{cases}
x_t^T Q x_t + u_t^T R u_t, & \|u_t\|_\infty \leq u_\mathrm{max}\\
+\infty & \text{otherwise.}
\end{cases}
\]
Unlike the LQR problem, in general, there is no known exact solution to the box-constrained
problem, analytical or otherwise.
Sophisticated methods can be used, however, to compute a lower bound
on the true optimal cost~\citep{wang2009performance}.

\paragraph{Policy.} Our COCP is
an ADP policy \eqref{e-adp} with a quadratic value function:
\begin{equation} \label{eq:box_lqr_policy}
\begin{array}{llll}
\phi(x) &=& \underset{u}{\argmin} & u^TRu + \|\theta(Ax+Bu)\|_2^2\\
&&\mbox{subject to} & \|u \|_\infty \leq u_\mathrm{max},
\end{array}
\end{equation}
with parameter $\theta \in\reals^{n\times n}$.
The lower bound found in~\citep{wang2009performance} yields a policy
that has this same form, for a particular value of $\theta$.

\begin{figure}
  \centering
  \includegraphics[]{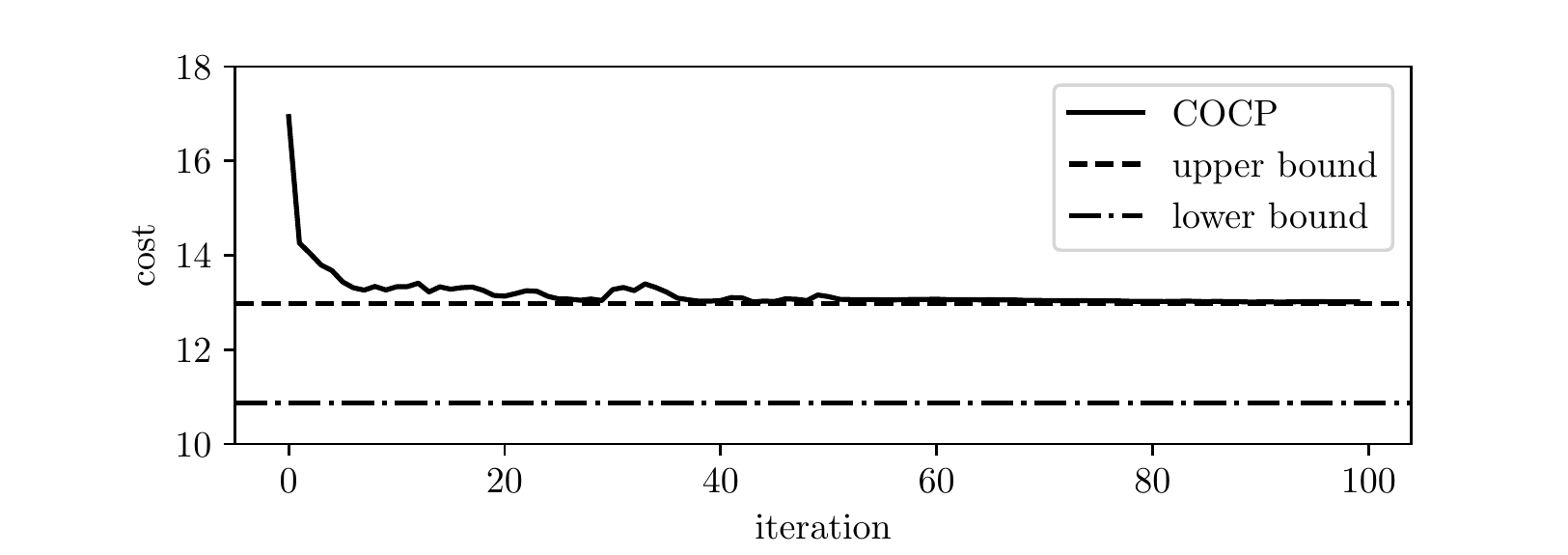}
  \caption{Tuning a box-constrained LQR policy.}
  \label{fig:lqr_constrained}
\end{figure}

\paragraph{Numerical example.}
We use $n=8$ states, $m=2$ inputs, $T=100$, $u_\mathrm{max}=0.1$, and
data generated as in the LQR example above.
The lower bounding technique from
\cite{wang2009performance} yields a lower bound on optimal
cost of around 11. It also suggests a particular value
of $\theta$, which gives average cost around 13,
an upper bound on the optimal cost that we suspect is the
true optimal average cost.
We initialize our COCP with $\theta=P^{1/2}$,
where $P$ comes from the cost-to-go function for the unconstrained (LQR) problem.
Figure~\ref{fig:lqr_constrained} plots the expected cost of our
COCP, and the expected cost of the upper and lower bounds suggested by~\cite{wang2009performance}.
Our method converges to roughly the same cost as the upper bound.

\subsection{Tuning a Markowitz policy to maximize utility}
In 1952, Markowitz introduced an optimization-based method
for the allocation of
financial portfolios \citep{markowitz1952portfolio}, which trades off
risk (measured as return variance), and (expected) return.
While the original formulation involved only a quadratic objective
and linear equality constraints (very much like LQR), with the addition of
other constraints and terms, %
Markowitz's method becomes a sophisticated COCP
\citep{grinold2000active,boyd2017multi}.
The parameters are the data that appear in the convex problem solved
to determine the trades to execute in each time period.

In this example, we learn the parameters in a Markowitz
policy to maximize a utility on the realized returns. We will use notation
from \cite{boyd2017multi}, representing the state by $w_t$, the control by $z_t
$, and the disturbance by $r_t$.

The portfolio under consideration has $n$ assets. The dollar value of
the portfolio in period $t$ is denoted by $v_t$, which we
assume to be positive. Our holdings in period $t$, normalized by the portfolio
value, are denoted by $w_t \in \reals^n$; the normalization ensures that
$\ones^Tw_t = 1$. The number $v_t(w_t)_i$ is the dollar value of our position
in asset $i$ ($(w_t)_i < 0$ corresponds to a short position).
In each period, we re-allocate our holdings by executing trades $z_t \in
\reals^n$, which are also normalized by $v_t$. Selling or
shorting asset $i$ corresponds to $(z_t)_i < 0$, and purchasing it corresponds
to $(z_t)_i
> 0$. Trades incur transaction costs $\kappa^T |z_t|$,
where $\kappa \in \reals_{++}^n$ (the set of positive $n$-vectors)
is the vector of transaction cost rates and
the absolute value is applied elementwise.
Shorting also incurs a cost, which we express by $\nu^T (w_t + z_t)_{-}$, where
$\nu \in \reals^n_{++}$ is the vector of stock loan rates and
$(\cdot)_{-}$ is the negative part. We impose the
condition that trades are self-financing, \ie, we must withdraw enough cash to
pay the transaction and shorting costs incurred by our trades. This can be
expressed as $\ones^T z_t + \kappa^T|z_t| + \nu^T (w_t + z_t)_{-} \leq 0$.

The holdings evolve according to the dynamics
\[
w_{t+1} = r_t \circ (w_t + z_t) / r_t^T (w_t + z_t)
\]
where $r_t \in \reals_{+}^n$ are the total returns (which are IID) and $\circ$ is
the elementwise product. The denominator in this expression is the return
realized by executing the trade $z_t$.

Our goal is to minimize the average negative utility of the realized returns, as
measured by a utility function $U : \reals \to \reals$. Letting $W$, $Z$ and $R$
denote the state, input, and disturbance trajectories, the cost function is
\[
\psi(W, Z, R) = \frac{1}{T+1}\sum_{t=0}^{T} -U(r_t^T (w_t + z_t)) + I(z_t),
\]
where $I : \reals^n \to \reals \cup \{+\infty\}$ enforces the
self-financing condition: $I(z_t)$ is $0$ when $\ones^T z_t + \kappa^T|z_t|
+ \nu^T (w_t + z_t)_{-} \leq 0$ and $+\infty$ otherwise.

\paragraph{Policy.}
We consider policies that compute $z_t$ as
\begin{equation*}
\begin{array}{llll}\label{marko}
\phi(w_t) &=& \underset{z}{\argmax} & \mu^T w^+ - \gamma\|S w^+\|_2^2  \\
&& \mbox{subject to} & w^+ = w_t + z \\
&&& \ones^T z + \kappa^T|z| + \nu^T (w^+)_{-} \leq 0,
\end{array}
\end{equation*}
with variables $w^+$ and $z$ and parameters $\theta = (\mu, \gamma, S)$, where
$\mu \in \reals^n$, $\gamma \in \reals_{+}$, and $S \in \reals^{n
\times n}$. In a Markowitz formulation, $\mu$ is set to the empirical mean
$\mu^{\mathrm{mark}}$ of the returns, and $S$ is set to
the square root of the return covariance $\Sigma^{\mathrm{mark}}$. With these
values for the parameters, the linear term in the objective represents the
expected return of the post-trade portfolio $w^+$, and the quadratic term
represents the risk. A trade-off between the risk and return is determined by
the choice of the risk-aversion parameter $\gamma$. We mention that it is
conventional to parametrize a Markowitz policy with a matrix $\Sigma \in
\mathbf{S}^n_{+}$, rewriting the quadratic term as ${w^+}^T \Sigma w^+$; as in
the LQR example, our policy is over-parametrized.

In addition to the self-financing condition, there are many other constraints one may
want to impose on the trade vector and the post-trade portfolio, including
constraints on the portfolio leverage and turnover, many of which are convex.
For various examples of such constraints, see~\citep[\S4.4, \S4.5]{boyd2017multi}.

\paragraph{Numerical example.}
\begin{figure}
\centering
\includegraphics{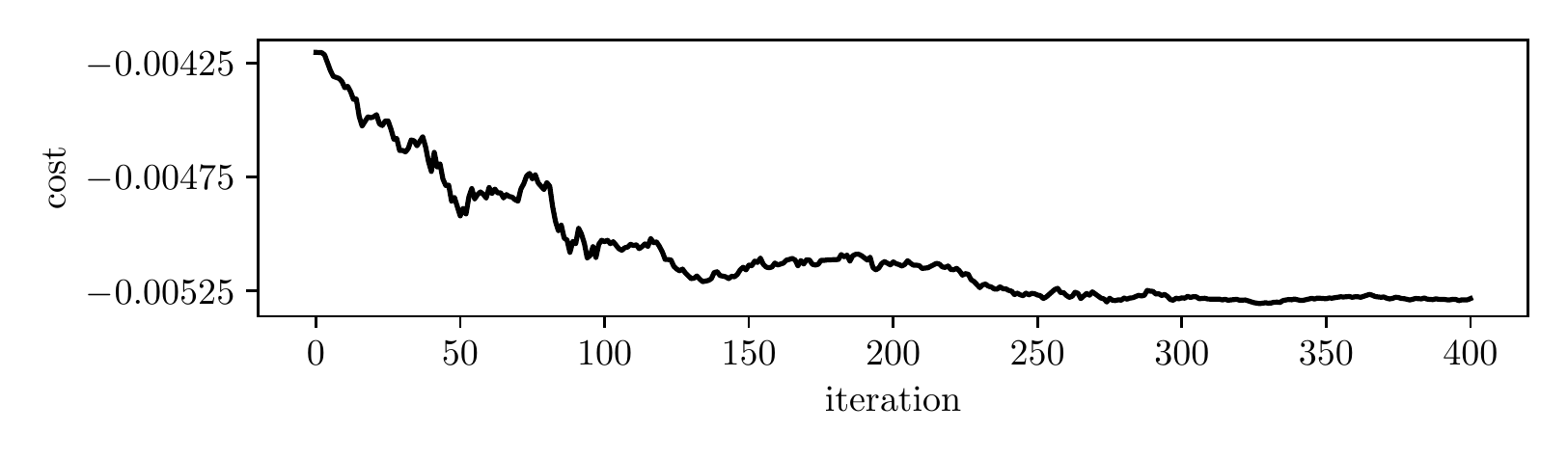}
\caption{Tuning a Markowitz policy.}
\label{fig-mark-train}
\end{figure}

We use $n=12$ ETFs as the universe of assets,
\begin{center}
AGG, VTI, VNQ, XLF, XLV, XLY, XLP, XLU, XLI, XLE, IBB, and ITA.
\end{center}
For the transaction rates and stock loan rates,
we use $\kappa = \nu = (0.001)\ones$, or $0.1$ percent. We assume the investor is somewhat
risk-averse, with utility function
\[
U(r) = \min(2(r-1), r-1).
\]

The policy is initialized with $\mu = \mu^{\mathrm{mark}}$, $S =
(\Sigma^{\mathrm{mark}})^{1/2}$, and $\gamma = 15$. Each simulation starts with the
portfolio obtained by solving
\begin{equation*}
\begin{array}{ll}\label{marko-init}
\mbox{maximize} & \mu^T w - \gamma\|S w\|_2^2 - \nu^T (w)_{-} \\
\mbox{subject to} & \ones^T w = 1,
\end{array}
\end{equation*}
with variable $w \in \reals^n$. The
portfolio evolves according to returns sampled from a log-normal distribution.
This distribution was fit to monthly returns (including dividends)
from Dec. 2006 through Dec. 2018, retrieved from the Center for Research in
Security Prices \citep{crsp}.

\begin{figure}
\centering
\adjustbox{valign=t}{\begin{minipage}[b]{.90\textwidth}
    \includegraphics[width=\linewidth]{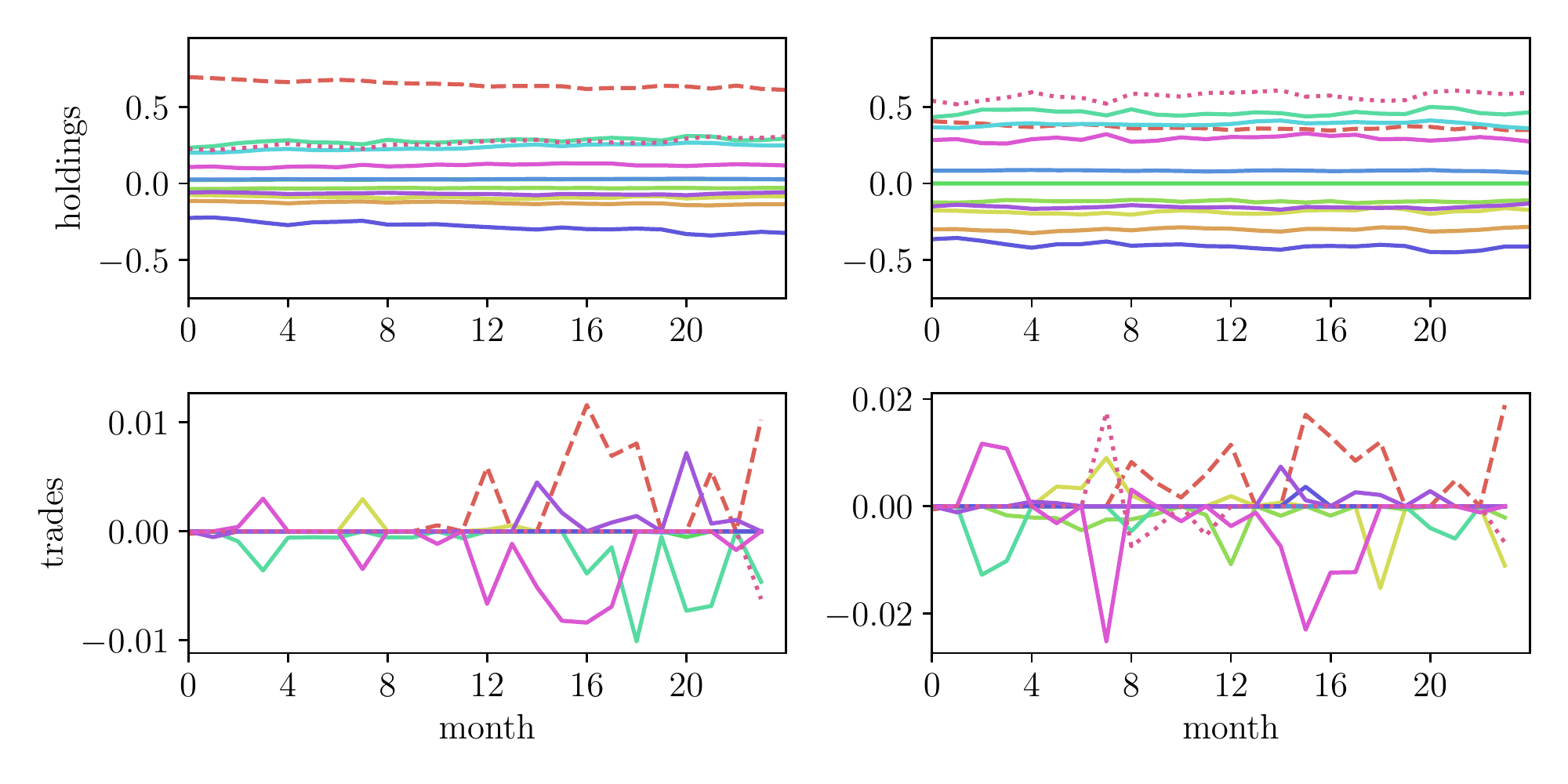}
\end{minipage}}
\adjustbox{valign=t}{\begin{minipage}{.09\textwidth}
    \vspace{31pt}
    \includegraphics[width=1.0\linewidth]{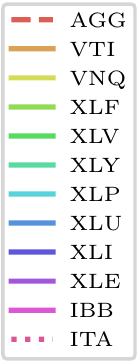}%
\end{minipage}}
\caption{Simulated holdings (top row) and trades (bottom row) for untuned (left column) and tuned (right column)
policies.}
\label{fig-mark}
\end{figure}

We train the policy using stochastic gradient descent over $400$ iterations,
with a horizon of $T = 24$ months and $K=10$ simulations to evaluate $\hat
J(\theta)$. (The step size is initialized to $10^{-3}$,
halved every $100$ iterations.) Figure~\ref{fig-mark-train} plots the
per-iteration cost on a held-out random seed
while training. The policy's
performance improved by approximately $32$ percent, decreasing from an initial
cost of $-0.004$ to $-0.0053$.

Figure~\ref{fig-mark} plots simulated holdings and trades before and after
tuning. Throughout the simulations, both the untuned and tuned policies
regulated or re-balanced their holdings to track the initial portfolio, making
small trades when their portfolios began to drift. The parameter $\mu$ was
adjusted from its initial value,
\[
(1.003, 1.006, 1.006, 1.002, 1.009, 1.009, 1.007, 1.006, 1.007, 1.004, 1.011, 1.011),
\]
to
\[
(0.999, 1.006, 1.005, 1.000, 1.001, 1.009, 1.008, 1.007, 1.009, 1.002, 1.014, 1.013).
\]
In particular, the entry corresponding to AGG, a bond ETF, decreased from $1.003$ to $0.999$,
and the entry for ITA, an aerospace and defense ETF, increased from $1.011$ to
$1.013$; this observation is consistent with the plotted simulated holdings.

Tuning had essentially no effect on $\gamma$, which decreased from $15$ to
$14.99$. The difference between $\Sigma^\mathrm{mark}$ and $S^T S$, however,
was significant: the median absolute percentage deviation between the entries
of these two quantities was $2.6$ percent.

\subsection{Tuning a vehicle controller to track curved paths}

We consider a vehicle moving relative to a smooth path,
with state and input
\[
x_t=(e_t,\Delta\psi_t,v_t,v^\mathrm{des}_t,\kappa_t),
\quad u_t=(a_t,z_t).
\]
Here, at time period $t$, $e_t$ is the lateral path deviation (\si{\meter}),
$\Delta\psi_t$ is the heading deviation from the path (\si{\radian}),
$v_t$ is the velocity (\si{\meter/\second}),
$v^\mathrm{des}_t$ is the desired velocity (\si{\meter/\second}),
$\kappa_t$ is the current curvature (\ie, inverse radius) of the path (\si{1/\meter}),
$a_t$ is the acceleration (\si{\meter/\second^2}), and
$z_t\coloneqq\tan(\delta_t) - L\kappa_t$,
where $\delta_t$ is the wheel angle (\si{\radian}) and $L$ is the vehicle's wheelbase (\si{\meter}).

\paragraph{Dynamics.}
We consider kinematic bicycle model dynamics in path coordinates \cite{gerdes},
discretized at $h=\SI{0.2}{\second}$,
with random processes for $v_t^{\mathrm{des}}$ and $\kappa_t$,
of the form
\begin{eqnarray}
e_{t+1} &=& e_t + h v_t\sin(\Delta\psi_t) + w_1, \quad w_1\sim\mathcal N(0, .01),\nonumber\\
{\Delta\psi}_{t+1} &=& {\Delta\psi}_{t} + h v_t \left(\kappa_t + \frac{z_t}{L} - \frac{\kappa_t}{1-e_t \kappa_t}\cos(\Delta\psi_t)\right) + w_2, \quad w_2\sim\mathcal N(0, .0001),\nonumber\\
v_{t+1} &=& v_t + ha_t + w_3 \nonumber, \quad w_3\sim\mathcal N(0, .01),\\
v_{t+1}^\mathrm{des} &=& v_t^\mathrm{des} w_4 + w_5 (1 - w_4), \quad w_4 \sim \mathrm{Bernoulli}(0.98), \quad w_5\sim \mathcal U(3, 6), \nonumber \\
\kappa_{t+1} &=& \kappa_t w_6 + w_7 (1 - w_6), \quad w_6 \sim \mathrm{Bernoulli}(0.95),\quad w_7\sim\mathcal N(0, .01). \nonumber
\end{eqnarray}
The disturbances $w_1,w_2,w_3$ represent uncertainty in our model,
and $w_4,\ldots,w_7$ form the random process for the desired speed
and path.

\paragraph{Cost.}
Our goal is to travel the desired speed ($v_t\approx v_t^\mathrm{des}$),
while tracking the path ($e_t\approx 0$, $\Delta\psi\approx 0$)
and expending minimal control effort ($a_t\approx 0$, $z_t \approx 0$).
We consider the cost
\[
\psi(X,U,W) = \frac{1}{T+1}\sum_{t=0}^T (v_t - v_t^\mathrm{des})^2 + \lambda_1 e_t^2 +
\lambda_2 \Delta\psi_t^2 + \lambda_3 |a_t| + \lambda_4 z_t^2 + I(a_t, z_t, \kappa_t),
\]
for positive $\lambda_1,\ldots,\lambda_4$ (with proper units),
where
\[
I(a, z,\kappa) = \begin{cases}0 &
|a| \leq a_\mathrm{max}, |z+L\kappa| \leq \tan(\delta_\mathrm{max}), \\
+\infty & \text{otherwise},
\end{cases}
\]
for given maximum acceleration magnitude $a_\mathrm{max}$ (\si{\meter/\second^2})
and maximum wheel angle magnitude $\delta_\mathrm{max}$ (\si{\radian}).

\paragraph{Policy.} We consider a COCP that computes $(a_t, z_t)$ as
\begin{equation*}
\begin{array}{llll}
\phi(x_t) &=& \underset{a, z}{\argmin} & \lambda_3 |a| + \lambda_4 z^2 + \|Sy\|_2^2 + q^T y\\
&&\mbox{subject to} &
y = \begin{bmatrix}
e_t + hv_t\sin(\Delta\psi_t) \\
{\Delta\psi}_{t} + h v_t \left(\kappa_t + \frac{z}{L} - \frac{\kappa_t}{1-e_t \kappa_t}\cos(\Delta\psi_t)\right) \\
v_t + ha - (0.98)v^\mathrm{des}_t - (0.02)4.5\\
y_1 + h v_t \sin(y_2 - hv_t\frac{z}{L}) + \frac{h^2v_t^2}{L} z,
\end{bmatrix}
\\
&&& |a| \leq a_\mathrm{max} \\
&&& |z + L \kappa_t| \leq \tan(\delta_\mathrm{max}),
\end{array}
\end{equation*}
with parameters $\theta=(S,q)$, where $S\in\reals^{4 \times 4}$ and $q\in\reals^4$.
The additional variable $y\in\reals^4$ represents relevant portions of the next state,
since $y_1=e_{t+1}$, $y_2=\Delta\psi_{t+1}$, $y_3=v_{t+1}-\Expect[v^\mathrm{des}_{t+1}]$,
and $y_4\approx e_{t+2}$ (since it assumes $a_t=0$).
Therefore, this COCP is an ADP policy and the term $\|Sy\|_2^2 + q^T y$ can
be interpreted as the approximate value function.

\begin{figure}
  \centering
  \includegraphics[width=\linewidth]{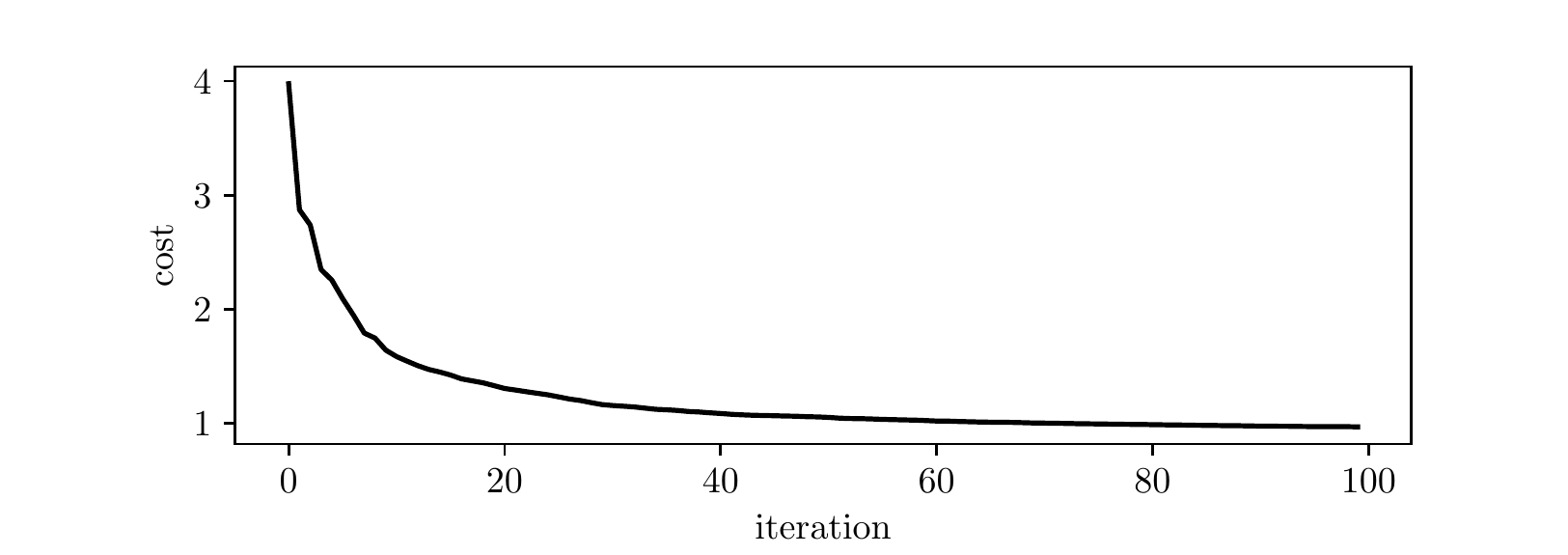}
\caption{Tuning a vehicle controller.}
\label{fig:vehicle_learning_curve}
\end{figure}

\begin{figure}
\begin{minipage}{.5\textwidth}
  \centering
  \includegraphics[]{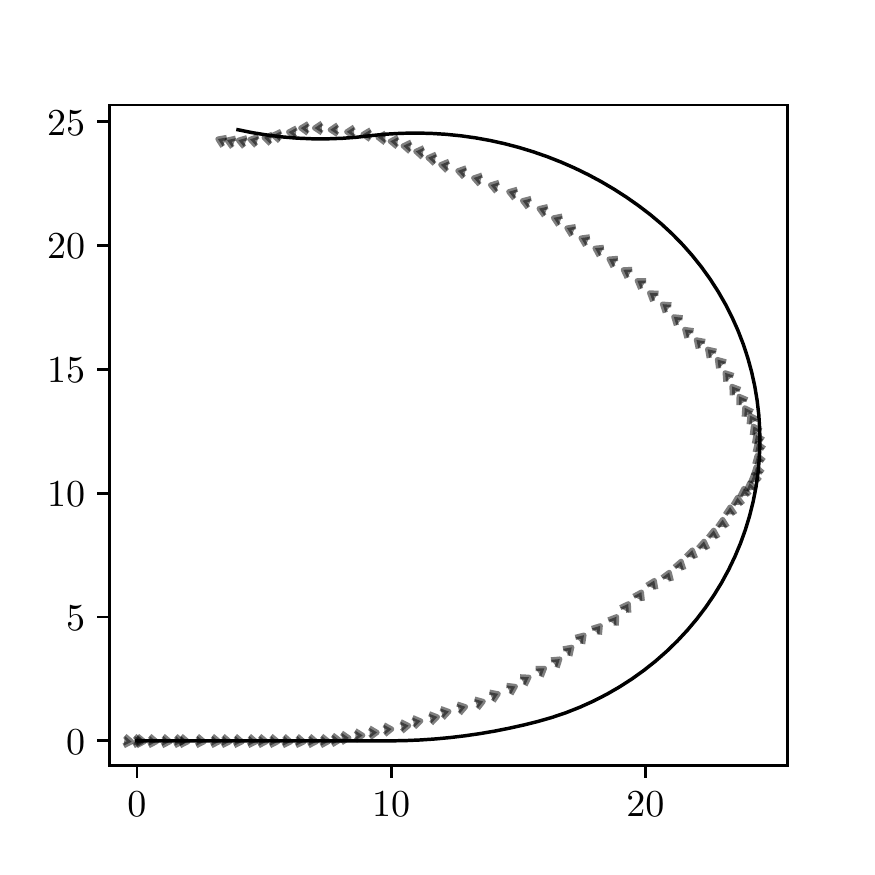}
\end{minipage}%
\begin{minipage}{.5\textwidth}
  \centering
  \includegraphics[]{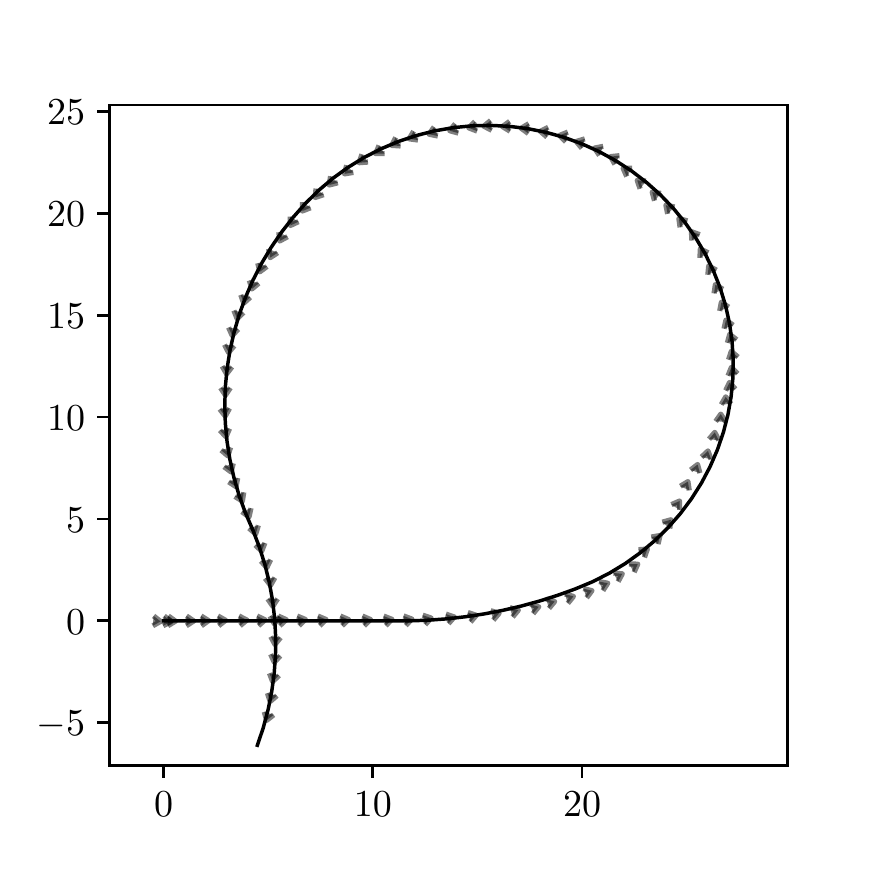}
\end{minipage}
\caption{Left: untuned policy. Right: tuned policy. Black line is the path and the gray
triangles represent the position and orientation of the vehicle. The tuned policy is able to track the path better and go faster.}
\label{fig:vehicle}
\end{figure}

\paragraph{Numerical example.}
We consider a numerical example with
\[
L=\SI{2.8}{\meter}, \quad \lambda_1=\lambda_2=1, \quad \lambda_3=\lambda_4=10, \quad a_\mathrm{max}=\SI{2}{\meter/\second^2}, \quad \delta_\mathrm{max}=\SI{0.6}{\radian},\quad T=100.
\]
We use the initial state $x_0=(.5, .1, 3, 4.5, 0)$.
We run the stochastic gradient method for $100$
iterations using $K=6$ simulations and a step size of 0.1.
We initialize the parameters with $S = I$ and $q = 0$.
Over the course of learning, the cost decreased from
$3.978$ to $0.971$.
Figure~\ref{fig:vehicle_learning_curve} plots
per-iteration cost on a held-out random seed
while training.
Figure~\ref{fig:vehicle} plots untuned and tuned sample paths
on a single held-out instance.
The resulting parameters are
\begin{equation*}
	S^TS = \begin{bmatrix}
  1.12 & 1.17 & -0.75 & 0.85\\
  1.17 & 3.82 & 0.46 & 3.13\\
  -0.75 & 0.46 & 13.07 & -0.29\\
  0.85 & 3.13 & -0.29 & 3.96\\
\end{bmatrix},
	\quad
q =
(-0, -0.04, -0.25, -0.04).
\end{equation*}

\subsection{Tuning a supply chain policy to maximize profit}
Supply chain management considers how to ship goods across a network of warehouses to maximize profit.
In this example, we consider a single-good supply chain with $n$ nodes representing interconnected warehouses linked to suppliers and consumers by $m$ directed links over which goods can flow.
There are $k$ links connecting suppliers to warehouses and $c$ links connecting warehouses to consumers.  The remaining $m - k - c$ links are internode links.

We represent the amount of good held at each node as $h_t \in \reals_+^n$
(the set of nonnegative $n$-vectors).
The prices at which we can buy the good from the suppliers
are denoted $p_t \in \reals^k_{+}$,
the (fixed) prices at which we can sell the goods to consumers
are denoted $r \in \reals^c_{+}$,
and the customer demand is denoted $d_t \in \reals^c_{+}$.
Our inputs are $b_t \in \reals^k_{+}$, the quantity of the good that we buy from the suppliers, $s_t \in \reals^c_{+}$, the quantity that we sell to the consumers, and $z_t \in \reals^{m -k -c}_{+}$,
the quantity that we ship across the internode links.
The state and inputs are
\begin{equation*}
	x_t = (h_t, p_t, d_t),\quad u_t=(b_t, s_t, z_t).
\end{equation*}
The system dynamics are
\begin{equation*}
	h_{t+1}= h_t + (A^{\mathrm{in}} - A^{\mathrm{out}})u_t,
\end{equation*}
where $A^{\mathrm{in}}\in\reals^{n \times m}$ and $A^{\mathrm{out}}\in\reals^{n \times m}$;
$A^{\mathrm{in(out)}}_{ij}$ is $1$ if link $j$ enters (exits) node $i$ and $0$ otherwise.

The input and state are constrained in several ways.
Warehouses have maximum capacities given by $h_{\rm max}\in\reals_+^n$, \ie, $h_t \leq h_{\rm max}$
(where the inequalities are elementwise),
and links have maximum capacities given by $u_{\rm max}\in\reals_+^m$, \ie, $u_t\leq u_{\rm max}$.
In addition, the amount of goods shipped out of a node cannot be more than the amount
on hand, or  $A^{\rm out} u_t \le h_t$.
Finally, we require that we sell no more than the demand, or $s_t \le d_t$.

We model the unknown future supplier prices and demands as random disturbances
$w_t = (p_{t+1}, d_{t+1})$ with joint log-normal distribution, \ie, $\log w_t =
(\log p_{t+1}, \log d_{t+1}) \sim \mathcal{N}(\mu, \Sigma)$.

The goal of our supply chain is to maximize profit,
which depends on several quantities.
Our payment to the suppliers is $p_t^T b_t$,
we obtain revenues $r^T s_t$ for selling the good to consumers, and
we incur a shipment cost $\tau^T z_t$, where $\tau \in \reals^{m-k-c}_{+}$
is the cost of shipping a unit of good across the internode links.
We also incur a cost for holding or storing $h_t$ in the warehouses;
this is represented by a quadratic function
$\alpha^T h_t + \beta^T h_t^2$, where $\alpha,\beta \in \reals^{n}_{++}$ and the square is elementwise.
Our cost is our average negative profit, or
\begin{equation*}
	\psi(X,U,W) = \frac{1}{T}\sum_{t=0}^{T-1} p_t^Tb_t -r^T s_t + \tau^T z_t + \alpha^T h_t + \beta^T h_t^2 + I(x_t, u_t).
\end{equation*}
Here, $I(x_t, u_t)$ enforces the constraints mentioned above; $I(x_t, u_t)$ is 0 if $x_t$ and $u_t$ lie in the set
\[
\left\{0 \le h_t \le h_{\rm max},
\quad 0 \le u_t \le u_{\rm max},
\quad A^{\rm out}u_t \le h_t,
\quad s_t \leq d_t
\right\},
\]
and $+\infty$ otherwise.

\paragraph{Policy.}
The policy seeks to maximize profit by computing $(b_t, s_t, z_t)$ as
\begin{equation*}
\begin{array}{llll}
	\phi(h_t, p_t, d_t) &=& \underset{b, s, z}{\argmax} & -p_t^Tb +r^T s - \tau^T z - \|S h^{+}\|_2^2 - q^Th^{+}\\
&&\mbox{subject to}
& h^{+} = h_t + (A^{\mathrm{in}}- A^{\mathrm{out}})(b, s, z) \\&&& 0 \le h^{+} \le h_{\rm max},\quad 0 \le (b, s, z) \le u_{\rm max},\\
&&& A^{\rm out}(b, s, z) \le h_t,\quad s \le d_t.
\end{array}
\end{equation*}
where the parameters are $\theta = (S, q)$ with $S\in \reals^{n \times n}$ and $q \in \reals^n$.
This COCP is an ADP policy and we can interpret the term $-\|S h^{+}\|_2^2 - q^Th^{+}$ as our approximate value function
applied to the next state.

\paragraph{Numerical example.}
\begin{figure}
\centering
\includegraphics[width=0.9\linewidth]{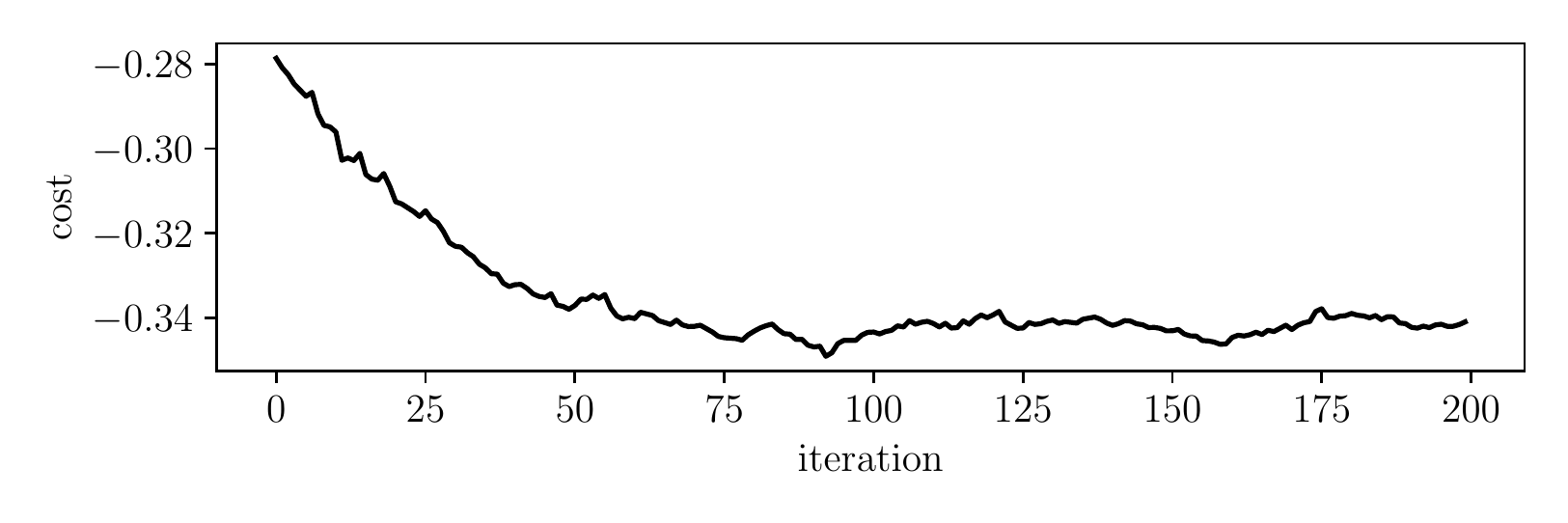}
\caption{Tuning a supply chain policy.}
\label{fig-supply-chain-train}
\end{figure}
We consider a supply chain over horizon $T=20$ with $n = 4$ nodes,
$m = 8$ links, $k = 2$ supply links, and $c = 2$ consumer links.
The initial value of the network storage is chosen uniformly between $0$ and $h_{\rm max}$, \ie, $h_0 \sim \mathcal{U}(0, h_{\rm max})$.
The log supplier prices and consumer demands have mean and covariance
\begin{equation*}
	\mu = (0.0, 0.1, 0.0, 0.4), \quad \Sigma = 0.04 I.
\end{equation*}
Therefore, the supplier prices have mean $(1.02,1.13)$ and the consumer demands have mean $(1.02,1.52)$.
The consumer prices are $r=(1.4)\ones$.
We set the maximum nodes capacity to $h_{\rm max} = (3)\ones$ and links capacity to $u_{\rm max}=(2)\ones$.
The storage cost parameters are $\alpha=\beta=(0.01)\ones$.
Node $1$ is connected to the supplier with lower average price and node $4$ to the consumer with higher demand.

We initialize the parameters of our policy to $S=I$ and $q=-h_{\rm max}$.
In this way, the approximate value function is centered at $h_{\rm max}/2$ so that we try to keep the storage of each node at medium capacity.

We ran our method over $200$ iterations, with $K=10$ using the stochastic gradient method with step size $0.05$.
Figure~\ref{fig-supply-chain-train} shows the per-iteration cost on a held-out random seed while training.
Over the course of training the cost decreased by $22.35$ percent from  $-0.279$ to $-0.341$.
The resulting parameters are
\begin{equation*}
	S^TS =
\begin{bmatrix}
  0.64 & 0.30 & 0.02 & -0.06\\
  0.30 & 1.44 & 0.32 & 0.30\\
  0.02 & 0.32 & 1.14 & 0.06\\
  -0.06 & 0.30 & 0.06 & 1.01\\
\end{bmatrix},
\quad
q=
(-3.05, -2.92, -2.97, -2.99).
\end{equation*}
The diagonal of $S^TS$ shows that the learned policy especially penalizes storing goods in nodes connected to more expensive suppliers, \eg, node $2$, or to consumers with lower demand, \eg, node $3$.
\begin{figure}
\centering
\includegraphics[width=0.9\linewidth]{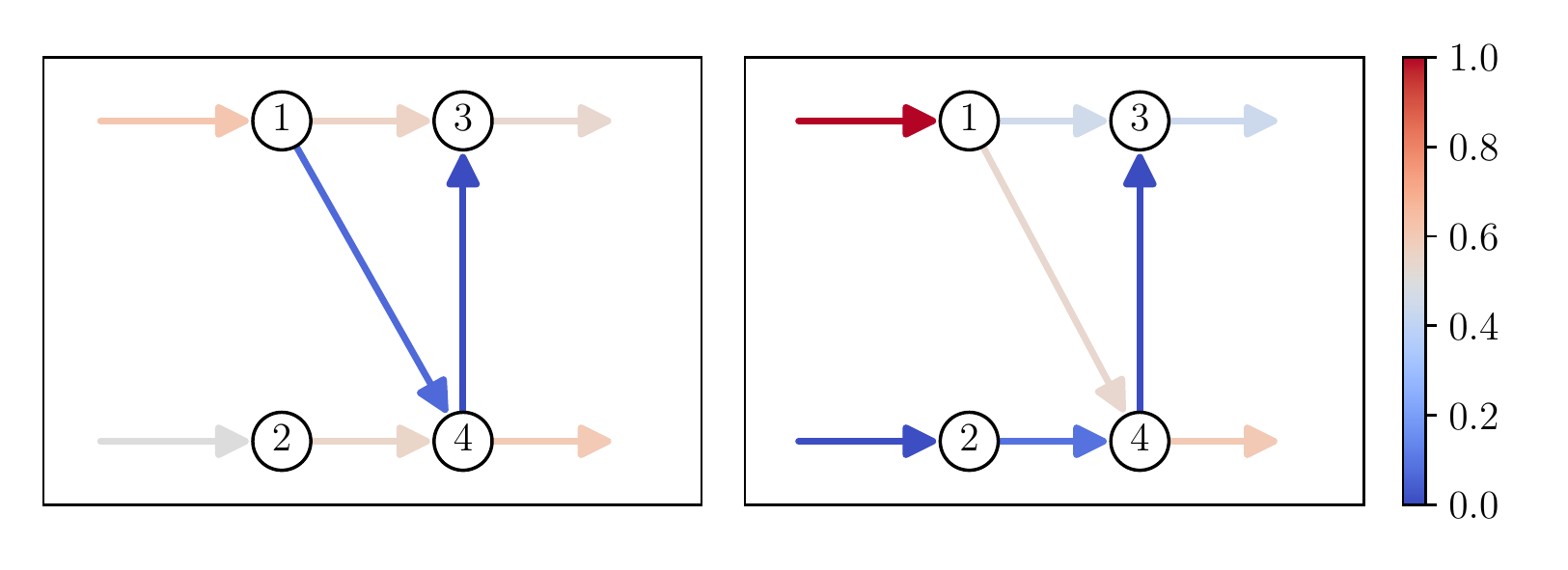}
\caption{Supply chain network. Left: untuned policy. Right: tuned policy. Colors indicate the normalized shipments between 0 and 1.}
\label{fig-supply-chain}
\end{figure}
\begin{figure}
\centering
\includegraphics[width=0.9\linewidth]{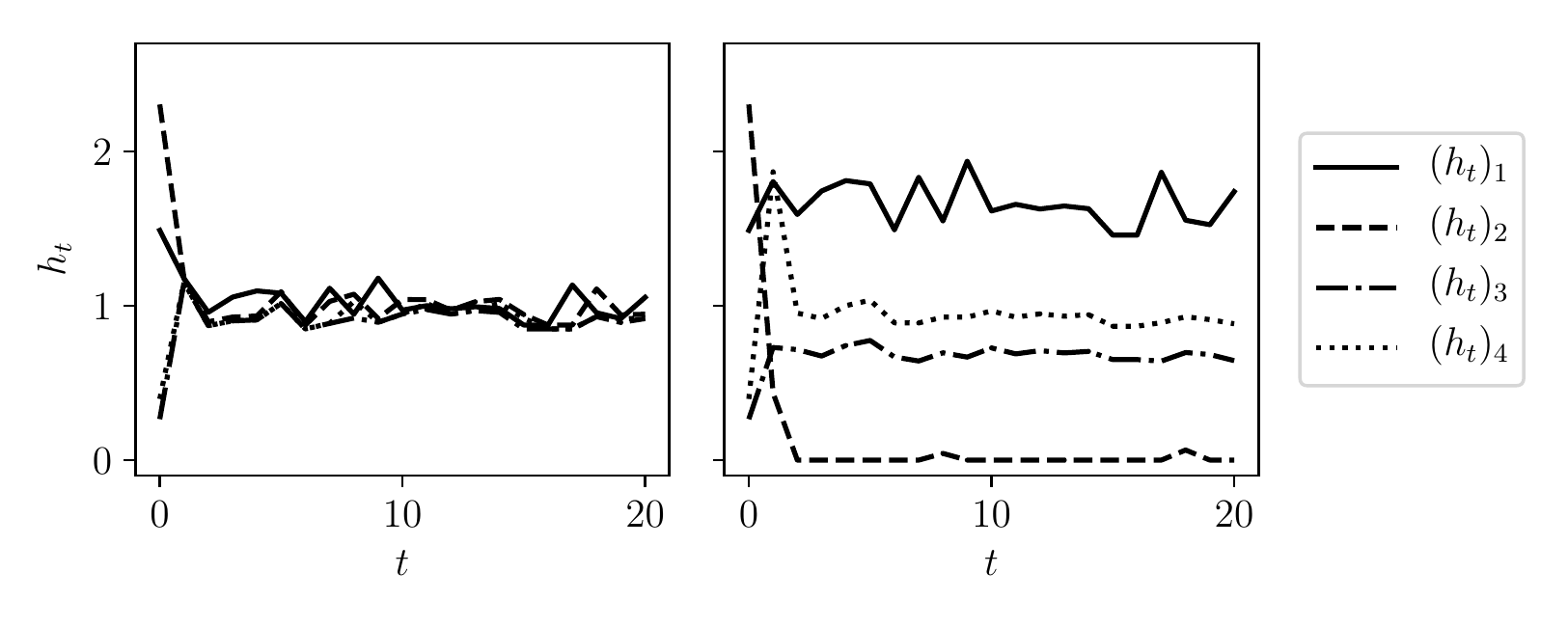}
\caption{Left: untuned policy. Right: tuned policy. Supply chain storage $h_t$ for each node over time.}
\label{fig-supply-chain-storage}
\end{figure}
Figure~\ref{fig-supply-chain} shows the supply chain structure and displays the
average shipment, normalized between $0$ and $1$;
figure~\ref{fig-supply-chain-storage} the simulated storage $h_t$ for the
untuned and tuned policy on a held-out random seed.

\section{Extensions and variations}%
\label{sec:extensions_and_variations}

\paragraph{Estimation.}
Our approach is not limited to tuning policies for control.
As we alluded to before, our approach can also be used
to learn convex optimization state estimators,
for example Kalman filters or moving horizon estimators.
The setup is exactly the same, in that we learn or tune
parameters that appear in the state estimation procedure
to maximize some performance metric.
(A similar approach was adopted in~\citep{barratt2019fitting},
where the authors fit parameters in a Kalman smoother to observed
data.)
Also, since COCPs are applied to the estimated state,
we can in fact jointly tune parameters in the COCP
along with the parameters in the state estimator.

\paragraph{Prediction.}
In an MPC policy, one could tune parameters in the function that predicts the
disturbances together with the controller parameters. As a specific example,
we mention that the parameters in a Markowitz policy, such as the expected
return, could be computed using a parametrized prediction function, and this
function could be tuned jointly with the other parameters in the COCP.

\paragraph{Nonconvex optimization control policies (NCOCPs).}
An NCOCP is an optimization-based control policy
that is evaluated by
solving a \emph{nonconvex} optimization problem.
Parameters in NCOCPs can be tuned in the same way that
we tune COCPs in this paper.
Although the solution to a nonconvex optimization problem might
be nonunique or hard to find,
one can differentiate a local solution map
to a smooth nonconvex optimization problem by implicitly
differentiating the KKT conditions~\citep{kuhn1951nonlinear}.
This is done in~\citep{amos2018mpc}, where the authors
define an MPC-based NCOCP.

\bibliography{refs}

\begin{thebibliography}{10}

\bibitem{abadi2016tensorflow}
M.~Abadi, P.~Barham, J.~Chen, Z.~Chen, A.~Davis, J.~Dean, M.~Devin,
  S.~Ghemawat, G.~Irving, M.~Isard, et~al.
\newblock {TensorFlow}: a system for large-scale machine learning.
\newblock In {\em OSDI}, volume~16, pages 265--283, 2016.

\bibitem{agrawal2019differentiable}
A.~Agrawal, B.~Amos, S.~Barratt, S.~Boyd, S.~Diamond, and J.~Z. Kolter.
\newblock Differentiable convex optimization layers.
\newblock In {\em Advances in Neural Information Processing Systems}, pages
  9558--9570, 2019.

\bibitem{diffcp2019}
A.~Agrawal, S.~Barratt, S.~Boyd, E.~Busseti, and W.~Moursi.
\newblock Differentiating through a cone program.
\newblock {\em Journal of Applied and Numerical Optimization}, 1(2):107--115,
  2019.

\bibitem{agrawal2018rewriting}
A.~Agrawal, R.~Verschueren, S.~Diamond, and S.~Boyd.
\newblock A rewriting system for convex optimization problems.
\newblock {\em Journal of Control and Decision}, 5(1):42--60, 2018.

\bibitem{amos2019differentiable}
B.~Amos.
\newblock {\em Differentiable optimization-based modeling for machine
  learning}.
\newblock PhD thesis, Carnegie Mellon University, 2019.

\bibitem{amos2018mpc}
B.~Amos, I.~Jimenez, J.~Sacks, B.~Boots, and J.~Z. Kolter.
\newblock Differentiable {MPC} for end-to-end planning and control.
\newblock In {\em Advances in Neural Information Processing Systems}, pages
  8299--8310, 2018.

\bibitem{amos2017input}
B.~Amos, L.~Xu, and J.~Z. Kolter.
\newblock Input convex neural networks.
\newblock In {\em Proc. Intl. Conf. on Machine Learning}, pages 146--155, 2017.

\bibitem{anderson1990}
B.~Anderson and J.~Moore.
\newblock {\em Optimal Control: Linear Quadratic Methods}.
\newblock Prentice-Hall, Inc., 1990.

\bibitem{anderson1953recent}
R.~Anderson.
\newblock Recent advances in finding best operating conditions.
\newblock {\em Journal of the American Statistical Association},
  48(264):789--798, 1953.

\bibitem{aastrom1993automatic}
K.~{\AA}str{\"o}m, T.~H{\"a}gglund, C.~Hang, and W.~Ho.
\newblock Automatic tuning and adaptation for pid controllers-a survey.
\newblock {\em Control Engineering Practice}, 1(4):699--714, 1993.

\bibitem{aastrom2013adaptive}
K.~{\AA}str{\"o}m and B.~Wittenmark.
\newblock {\em Adaptive Control}.
\newblock Courier Corporation, 2013.

\bibitem{osqp_codegen}
G.~Banjac, B.~Stellato, N.~Moehle, P.~Goulart, A.~Bemporad, and S.~Boyd.
\newblock Embedded code generation using the {OSQP} solver.
\newblock In {\em {IEEE} Conference on Decision and Control}, 2017.

\bibitem{barratt2018differentiability}
S.~Barratt.
\newblock On the differentiability of the solution to convex optimization
  problems.
\newblock {\em arXiv preprint arXiv:1804.05098}, 2018.

\bibitem{barratt2018stochastic}
S.~Barratt and S.~Boyd.
\newblock Stochastic control with affine dynamics and extended quadratic costs.
\newblock {\em arXiv preprint arXiv:1811.00168}, 2018.

\bibitem{barratt2019fitting}
S.~Barratt and S.~Boyd.
\newblock Fitting a kalman smoother to data.
\newblock {\em arXiv preprint arXiv:1910.08615}, 2019.

\bibitem{Bellman:1957}
R.~Bellman.
\newblock {\em Dynamic Programming}.
\newblock Princeton University Press, 1 edition, 1957.

\bibitem{bellman1957markovian}
R.~Bellman.
\newblock A markovian decision process.
\newblock {\em Journal of Mathematics and Mechanics}, pages 679--684, 1957.

\bibitem{bemporad2002}
A.~Bemporad, M.~Morari, V.~Dua, and E.~Pistikopoulos.
\newblock The explicit linear quadratic regulator for constrained systems.
\newblock {\em Automatica}, 38(1):3--20, 2002.

\bibitem{bengio1994learning}
Y.~Bengio, P.~Simard, and P.~Frasconi.
\newblock Learning long-term dependencies with gradient descent is difficult.
\newblock {\em IEEE Transactions on Neural Networks}, 5(2):157--166, 1994.

\bibitem{bergstra2012random}
J.~Bergstra and Y.~Bengio.
\newblock Random search for hyper-parameter optimization.
\newblock {\em Journal of Machine Learning Research}, 13(Feb):281--305, 2012.

\bibitem{bertsekas2017dynamic}
D.~Bertsekas.
\newblock {\em Dynamic Programming and Optimal Control}.
\newblock Athena Scientific, 4th edition, 2017.

\bibitem{bertsekas2019reinforcement}
D.~Bertsekas.
\newblock {\em Reinforcement Learning and Optimal Control}.
\newblock Athena Scientific, 1 edition, 2019.

\bibitem{bertsekas2004improved}
D.~Bertsekas, V.~Borkar, and A.~Nedi{\'c}.
\newblock Improved temporal difference methods with linear function
  approximation.
\newblock {\em Learning and Approximate Dynamic Programming}, pages 231--255,
  2004.

\bibitem{bertsekas1996neuro}
D.~Bertsekas and J.~Tsitsiklis.
\newblock {\em Neuro-dynamic Programming}, volume~5.
\newblock Athena Scientific, 1996.

\bibitem{spacex}
L.~Blackmore.
\newblock Autonomous precision landing of space rockets.
\newblock {\em The BRIDGE}, 26(4), 2016.

\bibitem{boot1963sensitivity}
J.~Boot.
\newblock On sensitivity analysis in convex quadratic programming problems.
\newblock {\em Operations Research}, 11(5):771--786, 1963.

\bibitem{borrelli2017predictive}
F.~Borrelli, A.~Bemporad, and M.~Morari.
\newblock {\em Predictive Control for Linear and Hybrid Systems}.
\newblock Cambridge University Press, 2017.

\bibitem{bottou2018optimization}
L.~Bottou, F.~Curtis, and J.~Nocedal.
\newblock Optimization methods for large-scale machine learning.
\newblock {\em {SIAM} Review}, 60(2):223--311, 2018.

\bibitem{boyd2017multi}
S.~Boyd, E.~Busseti, S.~Diamond, R.~Kahn, K.~Koh, P.~Nystrup, and J.~Speth.
\newblock Multi-period trading via convex optimization.
\newblock {\em Foundations and Trends{\textregistered} in Optimization},
  3(1):1--76, 2017.

\bibitem{boyd2004convex}
S.~Boyd and L.~Vandenberghe.
\newblock {\em Convex Optimization}.
\newblock Cambridge University Press, 2004.

\bibitem{busseti2018solution}
E.~Busseti, W.~Moursi, and S.~Boyd.
\newblock Solution refinement at regular points of conic problems.
\newblock {\em Computational Optimization and Applications}, 74:627--643, 2019.

\bibitem{crsp}
{Center for Research in Security Prices}.
\newblock Stock and security files, 2019.

\bibitem{chu2013code}
E.~Chu, N.~Parikh, A.~Domahidi, and S.~Boyd.
\newblock Code generation for embedded second-order cone programming.
\newblock In {\em European Control Conference}, pages 1547--1552. IEEE, 2013.

\bibitem{cochran2007sampling}
W.~Cochran.
\newblock {\em Sampling Techniques}.
\newblock John Wiley \& Sons, 2007.

\bibitem{conn2009introduction}
A.~Conn, K.~Scheinberg, and L.~Vicente.
\newblock {\em Introduction to Derivative-Free Optimization}, volume~8.
\newblock {SIAM}, 2009.

\bibitem{corless1988controller}
M.~Corless and G.~Leitmann.
\newblock Controller design for uncertain systems via {L}yapunov functions.
\newblock In {\em American Control Conference}, pages 2019--2025. IEEE, 1988.

\bibitem{cornuejols2006}
G.~Cornuejols and R.~T\"{u}t\"{u}nc\"{u}.
\newblock {\em Optimization Methods in Finance}.
\newblock Cambridge University Press, 2006.

\bibitem{de2003linear}
D.~De~Farias and B.~Van~Roy.
\newblock The linear programming approach to approximate dynamic programming.
\newblock {\em Operations Research}, 51(6):850--865, 2003.

\bibitem{diamond2016cvxpy}
S.~Diamond and S.~Boyd.
\newblock {CVXPY}: A {P}ython-embedded modeling language for convex
  optimization.
\newblock {\em Journal of Machine Learning Research}, 17(83):1--5, 2016.

\bibitem{domahidi2013ecos}
A.~Domahidi, E.~Chu, and S.~Boyd.
\newblock {ECOS}: An {SOCP} solver for embedded systems.
\newblock In {\em European Control Conference}, pages 3071--3076. IEEE, 2013.

\bibitem{gerdes}
C.~Gerdes.
\newblock {ME} 227 vehicle dynamics and control course notes, 2019.
\newblock Lectures 1 and 2.

\bibitem{goodfellow2016deep}
I.~Goodfellow, Y.~Bengio, and A.~Courville.
\newblock {\em Deep Learning}.
\newblock MIT Press, 2016.

\bibitem{gordon1995stable}
G.~Gordon.
\newblock Stable function approximation in dynamic programming.
\newblock In {\em Machine Learning}, pages 261--268. Elsevier, 1995.

\bibitem{grinold2000active}
R.~Grinold and R.~Kahn.
\newblock {\em Active Portfolio Management: A Quantitative Approach for
  Producing Superior Returns and Controlling Risk}.
\newblock McGraw Hill, 2000.

\bibitem{hansen2001completely}
N.~Hansen and A.~Ostermeier.
\newblock Completely derandomized self-adaptation in evolution strategies.
\newblock {\em Evolutionary computation}, 9(2):159--195, 2001.

\bibitem{jerez2014}
J.~Jerez, P.~Goulart, S.~Richter, G.~Constantinides, E.~C. Kerrigan, and
  M.~Morari.
\newblock Embedded online optimization for model predictive control at
  megahertz rates.
\newblock {\em IEEE Transactions on Automatic Control}, 59(12):3238--3251,
  2014.

\bibitem{kalman1960contributions}
R.~Kalman.
\newblock Contributions to the theory of optimal control.
\newblock {\em Boletin de la Sociedad Matematica Mexicana}, 5(2):102--119,
  1960.

\bibitem{keshavarz2012convex}
A.~Keshavarz.
\newblock {\em Convex methods for approximate dynamic programming}.
\newblock PhD thesis, Stanford University, 2012.

\bibitem{keshavarz2014quadratic}
A.~Keshavarz and S.~Boyd.
\newblock Quadratic approximate dynamic programming for input-affine systems.
\newblock {\em Intl. Journal of Robust and Nonlinear Control}, 24(3):432--449,
  2014.

\bibitem{kuhn1951nonlinear}
H.~Kuhn and A.~Tucker.
\newblock Nonlinear programming.
\newblock In {\em {B}erkeley {S}ymposium on {M}athematical {S}tatistics and
  {P}robability, 1950}, pages 481--492. University of California Press, 1951.

\bibitem{Kuindersma14}
S.~Kuindersma, F.~Permenter, and R.~Tedrake.
\newblock An efficiently solvable quadratic program for stabilizing dynamic
  locomotion.
\newblock In {\em Proc. Intl. on Robotics and Automation (ICRA)}, page
  2589{\textendash}2594, Hong Kong, China, 2014. IEEE, IEEE.

\bibitem{lillicrap2015continuous}
T.~Lillicrap, J.~Hunt, A.~Pritzel, N.~Heess, T.~Erez, Y.~Tassa, D.~Silver, and
  D.~Wierstra.
\newblock Continuous control with deep reinforcement learning.
\newblock {\em arXiv preprint arXiv:1509.02971}, 2015.

\bibitem{markowitz1952portfolio}
H.~Markowitz.
\newblock Portfolio selection.
\newblock {\em The Journal of Finance}, 7(1):77--91, 1952.

\bibitem{Mattingley:2012}
J.~Mattingley and S.~Boyd.
\newblock {CVXGEN}: A code generator for embedded convex optimization.
\newblock {\em Optimization and Engineering}, 13(1):1--27, 2012.

\bibitem{minorsky1922directional}
N.~Minorsky.
\newblock Directional stability of automatically steered bodies.
\newblock {\em Journal of the American Society for Naval Engineers},
  34(2):280--309, 1922.

\bibitem{movckus1975bayesian}
J.~Mo{\v{c}}kus.
\newblock On {B}ayesian methods for seeking the extremum.
\newblock In {\em Optimization Techniques IFIP Technical Conference}, pages
  400--404. Springer, 1975.

\bibitem{nedic2003least}
A.~Nedi{\'c} and D.~Bertsekas.
\newblock Least squares policy evaluation algorithms with linear function
  approximation.
\newblock {\em Discrete Event Dynamic Systems}, 13(1-2):79--110, 2003.

\bibitem{odonoghue2016conic}
B.~O'Donoghue, E.~Chu, N.~Parikh, and S.~Boyd.
\newblock Conic optimization via operator splitting and homogeneous self-dual
  embedding.
\newblock {\em Journal of Optimization Theory and Applications},
  169(3):1042--1068, 2016.

\bibitem{odonoghue2017scs}
B.~O'Donoghue, E.~Chu, N.~Parikh, and S.~Boyd.
\newblock {SCS}: {S}plitting conic solver, version 2.1.0.
\newblock \url{https://github.com/cvxgrp/scs}, 2017.

\bibitem{okada2017path}
M.~Okada, L.~Rigazio, and T.~Aoshima.
\newblock Path integral networks: End-to-end differentiable optimal control.
\newblock {\em arXiv preprint arXiv:1706.09597}, 2017.

\bibitem{pascanu2013difficulty}
R.~Pascanu, T.~Mikolov, and Y.~Bengio.
\newblock On the difficulty of training recurrent neural networks.
\newblock In {\em Proc. Intl. Conf. on Machine Learning}, pages 1310--1318,
  2013.

\bibitem{paszke2019pytorch}
A.~Paszke, S.~Gross, F.~Massa, A.~Lerer, J.~Bradbury, G.~Chanan, T.~Killeen,
  Z.~Lin, N.~Gimelshein, L.~Antiga, et~al.
\newblock {PyTorch}: An imperative style, high-performance deep learning
  library.
\newblock In {\em Advances in Neural Information Processing Systems}, pages
  8024--8035, 2019.

\bibitem{powell2007approximate}
W.~Powell.
\newblock {\em Approximate Dynamic Programming: Solving the Curses of
  Dimensionality}.
\newblock John Wiley \& Sons, 2007.

\bibitem{powell2012logistics}
W.~Powell, H.~Simao, and B.~Bouzaiene-Ayari.
\newblock Approximate dynamic programming in transportation and logistics: a
  unified framework.
\newblock {\em {EURO} Journal on Transportation and Logistics}, 1(3):237--284,
  2012.

\bibitem{rawlings2009model}
J.~Rawlings and D.~Mayne.
\newblock {\em Model Predictive Control: Theory and Design}.
\newblock Nob Hill Publishing, 2009.

\bibitem{robbins1951stochastic}
H.~Robbins and S.~Monro.
\newblock A stochastic approximation method.
\newblock {\em The Annals of Mathematical Statistics}, pages 400--407, 1951.

\bibitem{rosolia2017learning}
U.~Rosolia and F.~Borrelli.
\newblock Learning model predictive control for iterative tasks: {A}
  data-driven control framework.
\newblock {\em IEEE Transactions on Automatic Control}, 63(7):1883--1896, 2017.

\bibitem{rumelhart1988learning}
D.~Rumelhart, G.~Hinton, and R.~Williams.
\newblock Learning representations by back-propagating errors.
\newblock {\em Cognitive Modeling}, 5(3):1, 1988.

\bibitem{salimans2017evolution}
T.~Salimans, J.~Ho, X.~Chen, S.~Sidor, and I.~Sutskever.
\newblock Evolution strategies as a scalable alternative to reinforcement
  learning.
\newblock {\em arXiv preprint arXiv:1703.03864}, 2017.

\bibitem{schulman2017proximal}
J.~Schulman, F.~Wolski, P.~Dhariwal, A.~Radford, and O.~Klimov.
\newblock Proximal policy optimization algorithms.
\newblock {\em arXiv preprint arXiv:1707.06347}, 2017.

\bibitem{solis1981minimization}
F.~Solis and R.~Wets.
\newblock Minimization by random search techniques.
\newblock {\em Mathematics of Operations Research}, 6(1):19--30, 1981.

\bibitem{osqp}
B.~{Stellato}, G.~{Banjac}, P.~{Goulart}, A.~{Bemporad}, and S.~{Boyd}.
\newblock {OSQP}: {A}n operator splitting solver for quadratic programs.
\newblock {\em arXiv preprint arXiv:1711.08013}, 2019.

\bibitem{stellato2017}
B.~Stellato, T.~Geyer, and P.~Goulart.
\newblock High-speed finite control set model predictive control for power
  electronics.
\newblock {\em {IEEE} Transactions on Power Electronics}, 32(5):4007--4020,
  2017.

\bibitem{stewart2008}
G.~{Stewart} and F.~{Borrelli}.
\newblock A model predictive control framework for industrial turbodiesel
  engine control.
\newblock In {\em {IEEE} Conference on Decision and Control ({CDC})}, pages
  5704--5711, 2008.

\bibitem{sutton1988learning}
R.~Sutton.
\newblock Learning to predict by the methods of temporal differences.
\newblock {\em Machine Learning}, 3(1):9--44, 1988.

\bibitem{sutton2018reinforcement}
R.~Sutton and A.~Barto.
\newblock {\em Reinforcement Learning: An Introduction}.
\newblock MIT press, 2018.

\bibitem{tamar2017learning}
A.~Tamar, G.~Thomas, T.~Zhang, S.~Levine, and P.~Abbeel.
\newblock Learning from the hindsight plan--episodic {MPC} improvement.
\newblock In {\em IEEE Intl. on Robotics and Automation (ICRA)}, pages
  336--343, 2017.

\bibitem{wang2009performance}
Y.~Wang and S.~Boyd.
\newblock Performance bounds for linear stochastic control.
\newblock {\em Systems \& Control Letters}, 58(3):178--182, 2009.

\bibitem{wang2010fast}
Y.~Wang and S.~Boyd.
\newblock Fast evaluation of quadratic control-{L}yapunov policy.
\newblock {\em IEEE Transactions on Control Systems Technology},
  19(4):939--946, 2010.

\bibitem{wang2010}
Y.~{Wang} and S.~{Boyd}.
\newblock Fast model predictive control using online optimization.
\newblock {\em {IEEE} Transactions on Control Systems Technology},
  18(2):267--278, 2010.

\bibitem{wang2015approximate}
Y.~Wang, B.~O'Donoghue, and S.~Boyd.
\newblock Approximate dynamic programming via iterated {B}ellman inequalities.
\newblock {\em Intl. Journal of Robust and Nonlinear Control},
  25(10):1472--1496, 2014.

\bibitem{werbos1990backpropagation}
P.~Werbos.
\newblock Backpropagation through time: {W}hat it does and how to do it.
\newblock {\em Proc. IEEE}, 78(10):1550--1560, 1990.

\bibitem{williams1987reinforcement}
R.~Williams.
\newblock {\em Reinforcement-Learning Connectionist Systems}.
\newblock College of Computer Science, Northeastern University, 1987.

\end{thebibliography}

\end{document}